\documentclass[12pt]{article}
\usepackage[dvips]{epsfig}
\usepackage{graphicx}
\usepackage{amsthm,amsmath,amssymb}
\usepackage{mathtools}
\usepackage{latexsym}
\usepackage{graphics}
\usepackage{url}
\usepackage{placeins}
\usepackage{color}
\usepackage{bm}
\usepackage{comment}
\usepackage{statex}
\usepackage{pgfplots}

\usepackage[
            bookmarksnumbered=true,
            bookmarksopen=true,
            colorlinks=true,
            citecolor=green,
            linkcolor=green,
            anchorcolor=red,
            urlcolor=blue
            ]{hyperref}

\definecolor{blue}{rgb}{0,0,0.9}
\definecolor{red}{rgb}{0.9,0,0}
\definecolor{green}{rgb}{0,0.50,0.10}
\definecolor{violet}{rgb}{0.5804,0.0000,0.8275}

\def\@themcountersep{}

\makeatletter
\newcommand{\labeltext}[2]{%
  \@bsphack
  \csname phantomsection\endcsname 
  \def\@currentlabel{#1}{\label{#2}}%
  \@esphack
}
\makeatother

\newtheorem{THEO}{Theorem}[section]
\newtheorem{ALGo}[THEO]{Algorithm}
\newtheorem{CONJ}[THEO]{Conjecture}
\newtheorem{COND}[THEO]{Condition}
\newtheorem{CORO}[THEO]{Corollary}
\newtheorem{DEFI}[THEO]{Definition}
\newtheorem{EXAMP}[THEO]{Instance}
\newtheorem{INSTANCE}[THEO]{Instance}
\newtheorem{FACT}[THEO]{Fact}
\newtheorem{HYPO}[THEO]{Hypothesis}
\newtheorem{LEMM}[THEO]{Lemma}
\newtheorem{PROB}[THEO]{Problem}
\newtheorem{PROP}[THEO]{Proposition}
\newtheorem{REMA}[THEO]{Remark}
\newcommand{\theo}{\begin{THEO}}
\newcommand{\algo}{\begin{ALGo} \rm}
\newcommand{\cond}{\begin{COND}}
\newcommand{\conj}{\begin{CONJ}}
\newcommand{\coro}{\begin{CORO}}
\newcommand{\defi}{\begin{DEFI} \rm}
\newcommand{\examp}{\begin{EXAMP} \rm}
\newcommand{\instan}{\begin{INSTANCE} \rm}
\newcommand{\fact}{\begin{FACT}}
\newcommand{\hypo}{\begin{HYPO} \rm}
\newcommand{\lemm}{\begin{LEMM}}
\newcommand{\prob}{\begin{PROB} \rm}
\newcommand{\prop}{\begin{PROP}}
\newcommand{\rema}{\begin{REMA} \rm}
\newcommand{\etheo}{\end{THEO}}
\newcommand{\ealgo}{\end{ALGo}}
\newcommand{\econd}{\end{COND}}
\newcommand{\econj}{\end{CONJ}}
\newcommand{\ecoro}{\end{CORO}}
\newcommand{\edefi}{\end{DEFI}}
\newcommand{\eexamp}{\end{EXAMP}}
\newcommand{\einstan}{\end{INSTANCE}}
\newcommand{\efact}{\end{FACT}}
\newcommand{\ehypo}{\end{HYPO}}
\newcommand{\elemm}{\end{LEMM}}
\newcommand{\eprob}{\end{PROB}}
\newcommand{\eprop}{\end{PROP}}
\newcommand{\erema}{\end{REMA}}

\def\0{\mbox{\bf 0}}
\def\1{\mbox{\bf 1}}
\def\2{\mbox{\bf 2}}
\def\3{\mbox{\bf 3}}
\def\4{\mbox{\bf 4}}
\def\5{\mbox{\bf 5}}
\def\6{\mbox{\bf 6}}
\def\7{\mbox{\bf 7}}
\def\8{\mbox{\bf 8}}
\def\9{\mbox{\bf 9}}
\def\a{\mbox{\boldmath $a$}}
\def\b{\mbox{\boldmath $b$}}

\def\p{\mbox{\boldmath $p$}}

\def\ss{\mbox{\boldmath $s$}}

\def\u{\mbox{\boldmath $u$}}
\def\v{\mbox{\boldmath $v$}}

\def\x{\mbox{\boldmath $x$}}
\def\y{\mbox{\boldmath $y$}}

\def\A{\mbox{\boldmath $A$}}
\def\B{\mbox{\boldmath $B$}}
\def\C{\mbox{\boldmath $C$}}

\def\E{\mbox{\boldmath $E$}}

\def\H{\mbox{\boldmath $H$}}
\def\I{\mbox{\boldmath $I$}}

\def\L{\mbox{\boldmath $L$}}
\def\M{\mbox{\boldmath $M$}}

\def\O{\mbox{\boldmath $O$}}
\def\P{\mbox{\boldmath $P$}}
\def\Q{\mbox{\boldmath $Q$}}
\def\R{\mbox{\boldmath $R$}}
\def\S{\mbox{\boldmath $S$}}

\def\X{\mbox{\boldmath $X$}}
\def\Y{\mbox{\boldmath $Y$}}

\def\AC{\mbox{$\cal A$}}
\def\BC{\mbox{$\cal B$}}
\def\CC{\mbox{$\cal C$}}

\def\FC{\mbox{$\cal F$}}

\def\s0{\mbox{\scriptsize \boldmath $0$}}

\def\wFC{\mbox{$\widehat{\FC}$}}

\def\coneJ{\mathbb{J}}

\def\Real{\mathbb{R}}

\def\SymMat{\mathbb{S}}

\def\Integer{\mathbb{Z}}

\def\bGamma{\mbox{\boldmath $\Gamma$}}

\def\disk{\E^{\rm d}}
\def\hyperbola{\E^{\rm h}}
\def\parabola{\E^{\rm p}}
\def\linear{\E^{\rm \ell}}
\def\Example{Instance }

\begin{document}

\title{
Constructing QCQP Instances 
Equivalent to Their SDP Relaxations
}

\author{
\normalsize
Masakazu Kojima\thanks{Department of Data Science for Business Innovation ({\tt kojima@is.titech.ac.jp}).} \and \normalsize
Naohiko Arima\thanks{
	({\tt nao$\_$arima@me.com}).} \and \normalsize
Sunyoung Kim\thanks{Department of Mathematics, Ewha W. University, 52 Ewhayeodae-gil, Sudaemoon-gu, Seoul 03760, Korea 
			({\tt skim@ewha.ac.kr}). 
			 The research was supported  by   NRF 2021-R1A2C1003810.} \and \normalsize
}

\date{\today}

\maketitle 


\begin{abstract}
\noindent
General quadratically constrained quadratic programs (QCQPs) 
are challenging to solve as they are  known to be NP-hard. 
A popular approach to approximating QCQP solutions is to use 
semidefinite programming (SDP) relaxations.
It is well-known that the optimal value $\eta$ of the 
SDP relaxation problem bounds the optimal value 
$\zeta$ of the QCQP from below, {\it i.e.}, $\eta \leq \zeta$. 
The two problems are considered equivalent if $\eta = \zeta$. 
In the recent paper by Arima, Kim and Kojima \href{https://arxiv.org/abs/2409.07213}{[arXiv:2409.07213]}, 
a class of QCQPs that are equivalent to their SDP relaxations 
are proposed with no condition imposed on the quadratic objective function, which 
can be chosen arbitrarily.
In this work, we explore the construction of 
various QCQP instances 
within this class to complement the results in \href{https://arxiv.org/abs/2409.07213}{[arXiv:2409.07213]}.
Specifically, we first  construct QCQP instances with two variables and then extend them to 
higher dimensions. We also discuss how to compute an optimal QCQP solution  from the SDP relaxation. 
\end{abstract}


\vspace{0.5cm}



\section{Introduction}

We study a quadratically constrained quadratic program (QCQP) that minimizes a 
quadratic objective function in multiple real variables
over a feasible region represented 
by quadratic inequalities in the variables. 
In general, QCQPs are NP-hard \cite{MURTY87}
and become increasingly difficult to solve as the number of variables grows. 
It is well-known that 
the optimal value $\zeta$ of a QCQP is bounded by the optimal value $\eta$ of its SDP 
 relaxation from below  \cite{SHOR1987,Shor1990}. 
This property has frequently been utilized 
 for (approximately) solving the QCQP since the SDP can be solved numerically. We say that the QCQP and 
 the SDP relaxation are equivalent if $\eta = \zeta$. 
In \cite{ARIMA2024},  conditions on the feasible region 
were presented to ensure the equivalence of 
a QCQP and its SDP relaxation.
We first illustrate their conditions using the following figure.

\noindent
\hspace{-5mm}
\begin{minipage}{110mm}
\mbox{ \ } \hspace{1mm}
The unshaded area, including its boundary, 
represents the feasible region where 
two variables $u_1$ and $u_2$ can vary without limitation,
while the green  ellipsoidal regions indicate the restricted zones labeled 
$1,\ldots,7$. 
The feasible region is enclosed by the ellipsoid 8 within which 7 restricted zones 
are positioned. We regard the area outside of  the ellipsoid 8 as a restricted zone. 
Given a quadratic function $q(u_1,u_2)$ in 2 variables $u_1,u_2$, say 
$q(u_1,u_2) = u_1^2  + 2u_1u_2 +3u_2^2 -  2u_1 -2u_2$, we consider the 
optimization problem of minimizing $q(u_1,u_2)$ over the feasible region. 
The essential properties of this optimization problem are: 
\end{minipage}
\begin{minipage}{145mm}
\includegraphics[height=80mm]{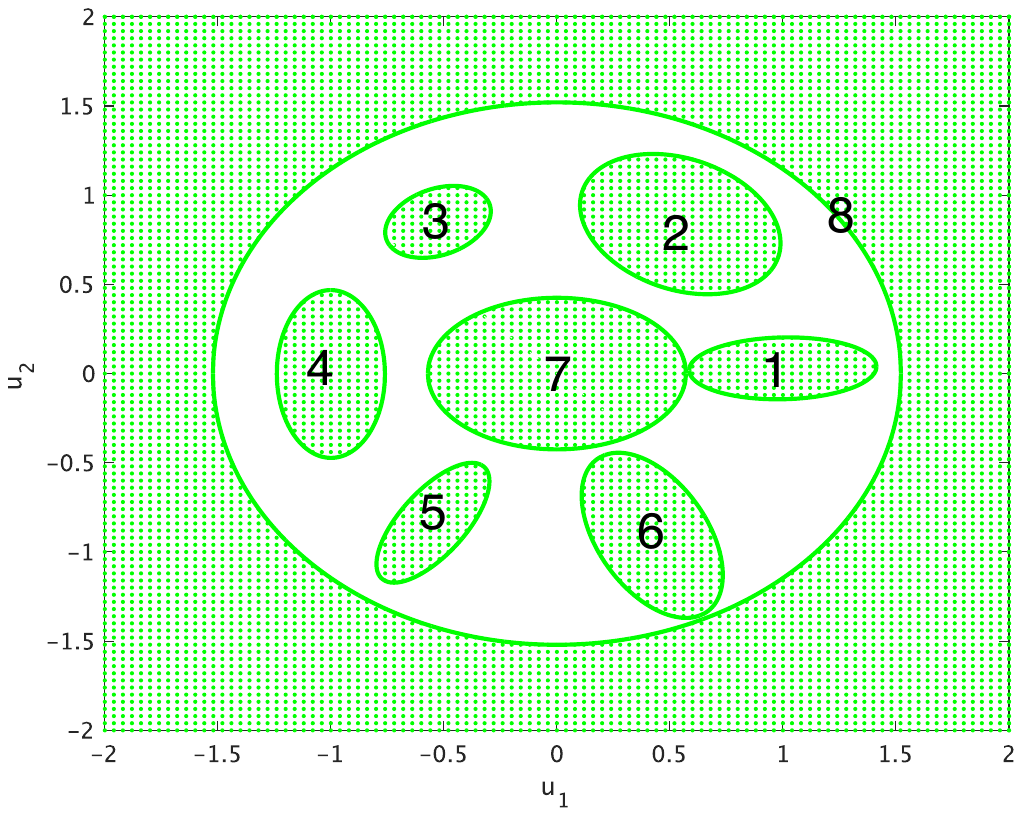}
\end{minipage}

\noindent
\begin{description}\vspace{-12mm}
\item{(a) } The objective function is a quadratic function in multiple variables. \vspace{-2mm}
\item{(b) } Each restricted zone is represented by a single quadratic inequality 
in multiple variables.\vspace{-2mm}
\item{(c) } Two distinct restricted zones could intersect with their boundaries but their interior never overlap 
(a more rigid condition is given as Condition (B)' below). \vspace{-2mm}
\end{description}

Properties (a) and (b) imply that this problem is a 
QCQP, while (c)  represents a condition  imposed on   the 
QCQPs that we consider throughout the paper. 
Since (c) imposes restrictions only on the constraints of QCQPs, the quadratic objective function can be chosen arbitrarily.
It was shown in \cite{ARIMA2024} that every QCQP satisfying (c) 
is equivalent to its SDP relaxation. 
The study in \cite{ARIMA2024}, however, focused mainly on the theoretical relationships conditions, 
including (c) that ensure the equivalence of QCQPs and their SDP relaxations, but 
did not fully address a variety of QCQP instances satisfying (c), nor
the construction of such instances. 
This paper aims to complement the work in \cite{ARIMA2024} by 
showing the  construction of  various QCQP instances satisfying (c). 
In addition to 
the instance mentioned above, 
Figures 1 through 13 
illustrate two-dimensional feasible regions satisfying condition (c). 
We also generate higher-dimensional QCQP instances
by extending the two-dimensional instances.

\smallskip

For the equivalence of a QCQP and its SDP relaxation, 
we  need to impose some conditions on the quadratic functions 
involved in the QCQP. 
Broadly speaking, there are three types of conditions 
that guarantee the equivalence between 
QCQPs and  their SDP relaxations.
The first type is  the convexity of quadratic objective 
function and quadratic inequality constraints. QCQPs  satisfying this condition are 
called as convex quadratic programs.  
The second type of conditions concerns the sign pattern of the coefficient matrices in 
both of quadratic objective and constraint functions. 
A class of QCQPs  satisfying this type of conditions was proposed in \cite{ZHANG2000} and 
has since been  extensively  studied in the literature  
\cite{AZUMA2022,KIM2003,SOJOUDI2014}. 

\smallskip

The focus of this paper is on the third type which imposes specific conditions only on 
quadratic inequality constraints, without any  conditions on the objective function.  
We note that condition (c) mentioned 
above falls under this type. 
If the constraints of a QCQP satisfy this type of conditions, then we can arbitrarily choose  
any quadratic objective function so that the QCQP is equivalent to its SDP relaxation. 
QCQPs satisfying this type of conditions were studied in \cite{ARGUE2023,
ARIMA2023,ARIMA2024,KIM2020}.  Specifically, we focus on 
Condition (B)' 
below  proposed in  \cite{ARIMA2024}.  
We note that Condition (B)' is translated into condition (c) above.

\medskip

For the subsequent discussion, we describe  the  standard form QCQP 
as the following: 
\begin{eqnarray}
\zeta & = & \inf \left\{ \x^T\Q\x : \x\in \Real^n, \ \x^T\B\x \geq 0 \ (\B \in \BC), \ \x^T\H\x = 1 \right\} \nonumber \\
& = & 
\inf \left\{ \Q\bullet \x\x^T : \x\in \Real^n, \ \B\bullet\x\x^T \geq 0 \ (\B \in \BC), \ \H\bullet\x\x^T = 1 \right\} 
\nonumber 
\\
& = & \inf \left\{ \Q\bullet \X : \X \in \bGamma^n, \ \B\bullet\X \geq 0 \ (\B \in \BC), \ \H\bullet\X = 1 \right\}. 
\label{eq:QCQP0} 
\end{eqnarray}
Here 
\begin{eqnarray*}
& & \Real^n : 
\mbox{the $n$-dimensional Euclidean space of column vectors $\x = (x_1,\ldots,x_n)$},\\
& & \SymMat^n : 
\mbox{the linear space of $n \times n$ symmetric matrices}, \\
& & \SymMat^n_+ \subseteq \SymMat^n : 
\mbox{the convex cone of $n\times n$ positive semidefinite matrices}, \\ 
& & \bGamma^n = \left\{\x\x^T \in \SymMat^n : \x \in \Real^n \right\}, \ 
 \Q \in \SymMat^n, \ \H \in \SymMat^n, \\
& & \BC : \mbox{a nonempty subset of $\SymMat^n$}, \\
& & \x^T : \mbox{the transposed row vector of $\x \in \Real^n$}, \\ 
& & \A \bullet \X = \sum_{i=1}^n\sum_{j=1}^n A_{ij}X_{ij}: \ 
 \mbox{the inner product of $\A, \ \X \in \SymMat^n$}.
\end{eqnarray*}
The set $\bGamma^n$ forms a cone in $\SymMat^n_+$; that is, for every 
 $\X \in \bGamma^n$ and $\lambda \geq 0$, it holds that $\lambda \X \in \bGamma^n$. 
 It is not convex unless $n = 1$. 
 We also know that  
$\SymMat^n_+ = \mbox{co}\bGamma^n$ (the convex hull of $\bGamma^n$). 
We may 
assume that $\BC$ is finite in this paper, though
Theorem~\ref{theorem:main} presented below remain valid  even when the cardinality of 
 $\BC$ is infinite.

\medskip

For every $\x \in \Real^n$, $\x\x^T \in \SymMat^n$ is 
an $n \times n$ rank-$1$ positive semidefinite 
matrix, and $\bGamma^n$ can be written as $\bGamma^n = \{ \X \in \SymMat^n: \mbox{rank}\X=1\}$. 
Hence we can rewrite \eqref{eq:QCQP0} as 
\begin{eqnarray*}
\zeta & = & \inf \left\{ \Q\bullet \X : \X \in \SymMat^n_+, \  \mbox{rank}\X=1, \ \B\bullet\X \geq 0 \ (\B \in \BC), \ \H\bullet\X = 1 \right\}. 
\end{eqnarray*}
If we remove rank$\X=1$ in the QCQP above (or relax $\bGamma^n$ by $\SymMat^n_+ \supseteq \bGamma^n$ in QCQP \eqref{eq:QCQP0}), we obtain an SDP relaxation of QCQP \eqref{eq:QCQP0}:  
\begin{eqnarray}
\eta & = & \inf \left\{ \Q\bullet \X : \X \in \SymMat^n_+, \ \B\bullet\X \geq 0 \ (\B \in \BC), \ \H\bullet\X = 1 \right\}. \label{eq:SDP0}
\end{eqnarray}
In general,  $\eta \leq \zeta$ holds. 

\medskip

As a special case of QCQP \eqref{eq:QCQP0}, we 
ocus on  the 
case where $\H \in \SymMat^n$ is the diagonal matrix 
with the diagonal entries $(0,\ldots,0,1) \in \Real^n$, {\it i.e.}, $\H = \mbox{diag}(0,\ldots,0,1) \in \SymMat^n$. 
Then the identity $\H\bullet\x\x^T = 1$ implies either $x_n = 1$ or $x_n = -1$. 
Since  
$\B\bullet\x\x^T = \B\bullet(-\x)(-\x)^T$ $(\B \in \BC)$ 
and $\Q\bullet\x\x^T = \Q\bullet(-\x)(-\x)^T$ for every $\x \in \Real^n$,  
we can fix $x_n = 1$ in this case. 
Therefore, we can rewrite QCQP \eqref{eq:QCQP0} with $\H = \mbox{diag}(0,\ldots,0,1) \in \SymMat^n$ 
as 
\begin{eqnarray}
\zeta & = & \inf \left\{ \begin{pmatrix}\u \\ 1\end{pmatrix}^T\Q\begin{pmatrix}\u \\ 1\end{pmatrix} : \u\in \Real^{n-1}, \ \begin{pmatrix}\u \\ 1\end{pmatrix}^T\B\begin{pmatrix}\u \\ 1\end{pmatrix} \geq 0 \ (\B \in \BC) \right\}. 
\label{eq:QCQP1} 
\end{eqnarray}

\medskip

To simplify the subsequent discussion, we introduce 
the following notation:
\begin{eqnarray*}
\coneJ_+(\B) \ \mbox{or } \coneJ_0(\B) & \equiv & 
\left\{ \X \in \SymMat^n_+ : \B\bullet \X \geq \ \mbox{or } = 0 
\right\}, \\ 
\coneJ_+(\BC) & \equiv & \bigcap_{\B \in \BC} \coneJ_+(\B)
 = \left\{ \X \in \SymMat^n_+ : \B\bullet \X \geq 0 \ 
(\B \in \BC) \right\}, \\
\B_{\geq}, \ \B_{\leq} \  \mbox{ or }  \B_{<} &\equiv& 
\left\{ \u \in \Real^{n-1} : \begin{pmatrix}\u \\ 1 \end{pmatrix}^T \B \begin{pmatrix}\u \\ 1 \end{pmatrix}
\geq, \ \leq  \ \mbox{or } < 0 \right\}, \\
\BC_\geq  & \equiv & \bigcap_{\B \in \BC} \B_\geq 
 =  
\left\{ \u \in \Real^{n-1} : \begin{pmatrix}\u \\ 1 \end{pmatrix}^T \B \begin{pmatrix}\u \\ 1 \end{pmatrix} \geq 0 
\ (\B \in \BC) \right\} \\ 
& & 
\mbox{(the feasible region of QCQP \eqref{eq:QCQP1})}
\end{eqnarray*}
for  every $\B \in \SymMat^n$ and $\BC \subseteq \SymMat^n$. 
Using the above notation, we can rewrite QCQP \eqref{eq:QCQP0}, its SDP relaxation \eqref{eq:SDP0} 
and QCQP \eqref{eq:QCQP1} as 
\begin{eqnarray}
\zeta & = & \inf \left\{ \Q\bullet \X : \X \in \bGamma^n \cap \coneJ_+(\BC), \ \H\bullet\X = 1 \right\}, 
\label{eq:QCQP2}
 \\ 
\eta & = & \inf \left\{ \Q\bullet \X : \X \in  \coneJ_+(\BC), \ \H\bullet\X = 1 \right\}, 
\label{eq:SDP2}
\\
\zeta & = & 
 \inf \left\{ \begin{pmatrix}\u \\ 1\end{pmatrix}^T\Q\begin{pmatrix}\u \\ 1\end{pmatrix} : 
 \u \in \BC_\geq \right\}, \label{eq:QCQP3}
\end{eqnarray}
respectively. 

\medskip

We now consider the following conditions on $\BC \subseteq \SymMat^n$. 

\begin{description} \vspace{-2mm} 
\item[(B) ] $\coneJ_0(\B) \subseteq \coneJ_+(\BC)$ $(\B \in \BC)$.  
\vspace{-2mm} 
\item[(B)' ]  
$\A_< \cap \B_\leq = \emptyset$ for every distinct pair of $\A\in\BC$ and $\B \in\BC$.\vspace{-2mm} 
\item[(C)' ]  $\B_< \not= \emptyset$ for every $\B \in \BC$.\vspace{-2mm}
\item[(D) ] There exist $\alpha_{B} > 0$ $(\B \in \BC)$ such that 
$\alpha_A\A + \alpha_B\B \in \SymMat^n_+$ for every distinct pair of $\A\in\BC$ and $\B \in\BC$.
\vspace{-2mm} 
\end{description}

Conditions (B), (B)' and (C)' were introduced in \cite[Section 1]{ARIMA2024}, 
and (D) is a special case of a sufficient condition discussed in \cite[Lemma 3.4]{ARIMA2024} 
for (B). 

\theo \label{theorem:main} \mbox{ \ }\vspace{-2mm}
\begin{description}
\item{(i) } 
Let $\Q, \ \H \in \SymMat^n$. 
Assume that 
\begin{eqnarray*}
\coneJ_+(\BC) \in 
 \wFC(\bGamma^n) \equiv  
 \left\{ 
 \mbox{
 closed convex cone $\coneJ \subseteq \SymMat^n_+ : 
\coneJ = \mbox{co}(\coneJ \cap \bGamma^n)$
} 
\right\},
\end{eqnarray*}
where $\mbox{co}(\coneJ \cap \bGamma^n)$ denotes the convex hull of 
$\coneJ \cap \bGamma^n$. Then 
\begin{eqnarray*}
-\infty < \eta \ \Leftrightarrow  -\infty < \eta = \zeta
\end{eqnarray*}
in QCQP \eqref{eq:QCQP2} and its SDP relaxation \eqref{eq:SDP2}. 
($\coneJ \in \wFC(\bGamma^n)$ is called {\em ROG (Rank-One Generated)} in the literature \cite{ARGUE2023,HILDEBRAND2016}).
\vspace{-2mm}
\item{(ii) } 
Condition (B) $\Rightarrow$ $\coneJ_+(\BC) \in  \wFC(\bGamma^n)$. 
\vspace{-2mm}
\item{(iii) }  Condition (B)' and (C)' $\Rightarrow$ $\coneJ_+(\BC) \in  \wFC(\bGamma^n)$.
\vspace{-2mm}
\item{(iv) }  Conditiion (D) $\Rightarrow$ Conditions (B) and (B)'.
 \vspace{-2mm}
\end{description}
\etheo
\proof{See \cite[Theorems 3.1]{KIM2020} for (i),  \cite[Theorems 1.2]{ARIMA2024} for (ii)
and \cite[Theorem 1.3]{ARIMA2024} for (iii). 
(iv) can be proved easily. To prove (D) $\Rightarrow$ (B), let $\B \in \BC$ and 
$\X \in \coneJ_0(\B)$. Then we see that  
\begin{eqnarray*}
0 \leq (\alpha_{A} \A + \alpha_{B}\B) \bullet \X =  \alpha_{A} \A\bullet\X  \ \mbox{for every }
\A \in \BC, 
\end{eqnarray*}
which implies $\X \in \coneJ_+(\BC)$. 
To prove (D) $\Rightarrow$ (B)', assume on the contrary that $\u \in \A_< \cap \B_\leq
 \not= \emptyset$ 
for some distinct $\A, \B \in \BC$. Then
\begin{eqnarray*}
0 > \begin{pmatrix} \u \\ 1 \end{pmatrix}^T \alpha_A \A \begin{pmatrix} \u \\ 1 \end{pmatrix}
+ \begin{pmatrix} \u \\ 1 \end{pmatrix}^T \alpha_B \B \begin{pmatrix} \u \\ 1 \end{pmatrix} 
= \begin{pmatrix} \u \\ 1 \end{pmatrix}^T (\alpha_A\A+ \alpha_B\B) \begin{pmatrix} \u \\ 1 \end{pmatrix} \geq 0, 
\end{eqnarray*}
which is a contradiction. \\
\qed
} 

\smallskip

\rema Assertion (iii) above can be strengthened to \vspace{-2mm}
\begin{description} 
\item{(iii)' } Condition (B)' $\Rightarrow$ $\coneJ_+(\BC) \in  \wFC(\bGamma^n)$.
\vspace{-2mm}
\end{description} 
In fact, we can prove that 
Condition (C)' is equivalent to $\BC \cap \SymMat^n_+ = \emptyset$. If $\B \in \SymMat^n_+$ 
then $\coneJ_+(\B) = \SymMat^n_+$ and $\B_\geq = \Real^{n-1}$. Hence 
$\coneJ_+(\BC\backslash\SymMat^n_+) = \coneJ_+(\BC)$ and 
$(\BC\backslash\SymMat^n_+)_\geq = \BC_\geq$.  
Therefore, we will not explicitly refer to Condition (C)’  in the following discussion, 
assuming that it is satisfied.  
\erema
It should be noted that $\coneJ_+(\BC) \in \wFC(\bGamma^n)$ does not 
 involve  $\Q\in \SymMat^n$ and 
$\H \in \SymMat^n$. By Theorem~\ref{theorem:main}, if $\BC$ satisfies on 
Condition (B)' 
imposed on QCQP \eqref{eq:QCQP3}, 
which is a special case of QCQP \eqref{eq:QCQP2} with 
$\H=\mbox{diag}(0,\ldots,0,1)\in \SymMat^n_+$, then $ -\infty < \eta \ \Leftrightarrow  -\infty < \eta = \zeta$ holds in QCQP \eqref{eq:QCQP2} and its SDP relaxation \eqref{eq:SDP2} 
for every $\Q, \ \H \in \SymMat^n$. 
In particular, 
if we take a positive definite matrix for $\H \in \SymMat^n$, 
the feasible regions of both QCQP \eqref{eq:QCQP2} and its SDP relaxation \eqref{eq:SDP2} 
are bounded; hence $-\infty < \eta$. In this case, we have $\eta = \zeta$. 
When the feasible region of SDP relaxation \eqref{eq:SDP2} is unbounded, however, 
$-\infty = \eta < \zeta$ 
may happen even if $\coneJ_+(\BC) \in  \wFC(\bGamma^n)$ holds. We will show such a case in \Example 2.7.

\subsection*{Main contribution of the paper}

Our first contribution is to present various QCQP instances satisfying Condition (B)'  by introducing a 
systematic method to construct such QCQPs. 
This work is based on the authors' previous work \cite{ARIMA2023,ARIMA2024,KIM2020}. In \cite{ARIMA2024}, 
Condition (B)' was proposed as a sufficient condition for the equivalence between QCQP and its SDP relaxation. 
A detailed theoretical analysis was also conducted to examine its relationship with several previously proposed sufficient conditions 
in  \cite{ARIMA2023,KIM2020}.   
While the paper \cite{ARIMA2024} provided
some examples of QCQPs satisfying Condition (B)', 
it remains unclear what kinds 
of problems the entire class of such QCQPs include. 
This work seeks to bridge this gap by presenting 
various QCQP instances within this class. Specifically, for two-dimensional QCQPs, 
we provide instances that fully utilize the geometric properties of Condition (B)’ to make them more intuitive and easier to understand. 
Based on these, it is demonstrated that various instances of high-dimensional QCQPs can be constructed freely. 
These findings deepen the understanding of Conditions (B)' proposed in 
\cite{ARIMA2024}  and are expected to facilitate the broader application of QCQPs that satisfy Conditions (B)'.

The second contribution of this work is the introduction of a numerical method to compute the optimal solution of a QCQP satisfying Condition (B)
from the optimal solution of
its SDP relaxation. The method is based on \cite[Lemma 2.2]{YE2003} and its constructive proof, and is
illustrated with numerical 
results to demonstrate its effectiveness.

\subsection*{Organization of the paper}

In Section 2, we deal with the case $n=3$. In this case, \eqref{eq:QCQP3} is a 
QCQP in the $2$-dimensional variable vector $\u$. After demonstrating how to construct 
 basic quadratic constraints in $\u$, we combine them to generate
 $7$ instances of $\BC$  satisfying Conditions (B)' and (D). 
 In Section 3, we discuss how to combine those instances for constructing 
 higher-dimensional $\BC$ satisfying  Condition (D). In Section~4, we show how to compute 
 an optimal solution $\widetilde{\X} \in \bGamma^n$ of QCQP~\eqref{eq:QCQP2} from 
 an optimal solution $\overline{\X} \in \SymMat^n_+$ of its SDP relaxation~\eqref{eq:SDP2} 
 under Condition (B) and some additional 
 assumption.


\section{QCQP \eqref{eq:QCQP3} with $\u \in \Real^2$}

Throughout this section, we assume that $n=3$ and $\u \in \Real^2$ in QCQP \eqref{eq:QCQP3}. 
Sections 2.1 and 2.2 describe quadratic constraints for Section 2.3  
where $7$ instances of $\BC$ satisfying Conditions (B)' and (D) are presented. 
We introduce $4$ types of basic quadratic constraints in Section 2.1, 
and $3$ types of linear transformations for them in Section 2.2. 

\medskip
 
\subsection{Basic quadratic constraints}

\medskip

\noindent
{\bf A disk constraint: }
\begin{eqnarray*}
& &\disk(r) \equiv \begin{pmatrix} 1 & 0 & 0 \\ 0 & 1 & 0 \\ 0 & 0 & -r^2 \end{pmatrix}, \
\disk(r)_{\geq} \ \mbox{or } \disk(r)_{\leq} 
= \left\{ \u \in \Real^2 : u_1^2 + u_2^2 -r^2 \geq \ \mbox{or } \leq 0 \right\},  
\end{eqnarray*}
where $r > 0$ denotes a parameter. See Figure 1. 

\medskip

\begin{figure}[t!]   \vspace{-30mm} 
\begin{center}
\includegraphics[height=110mm]{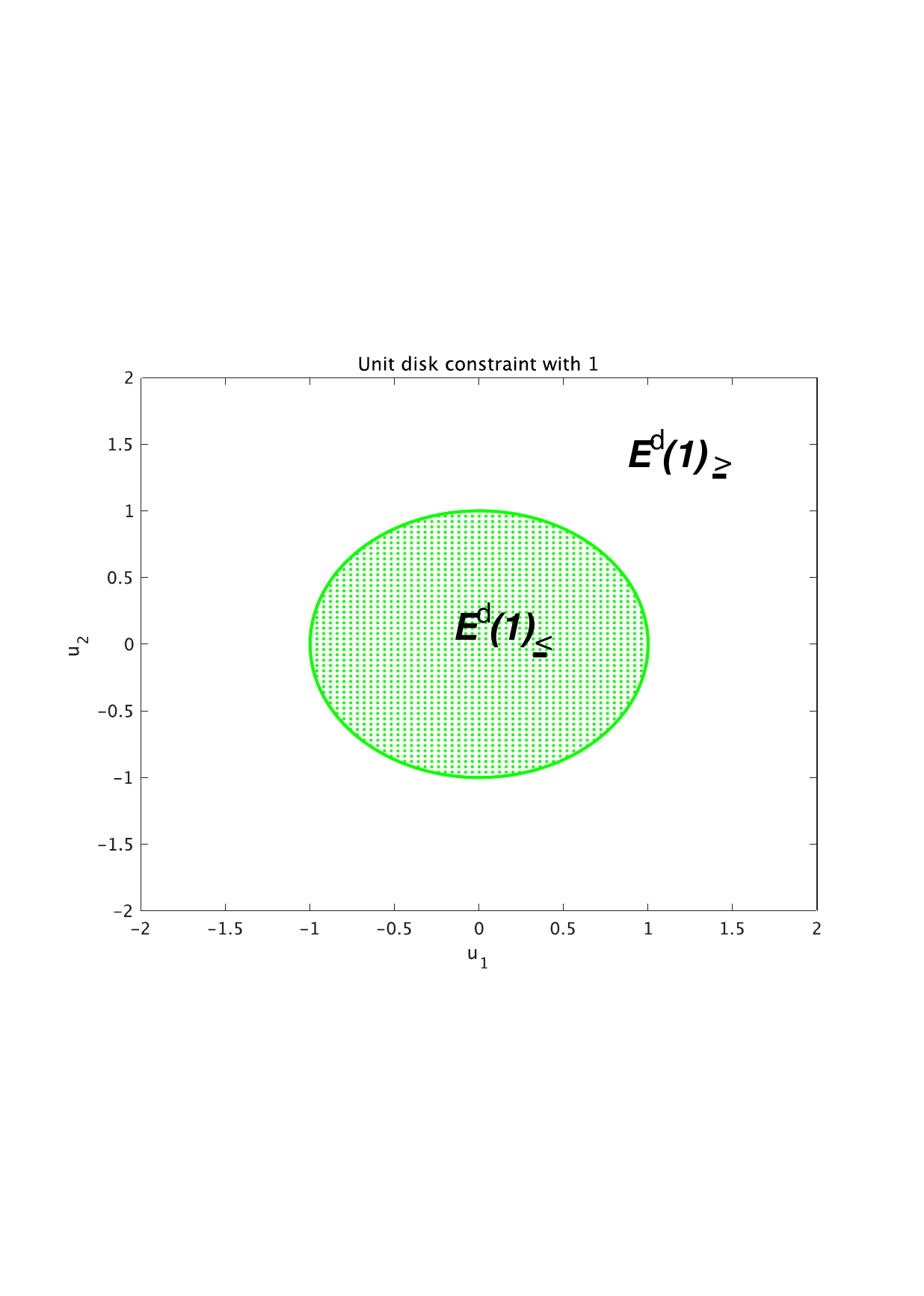}
\includegraphics[height=110mm]{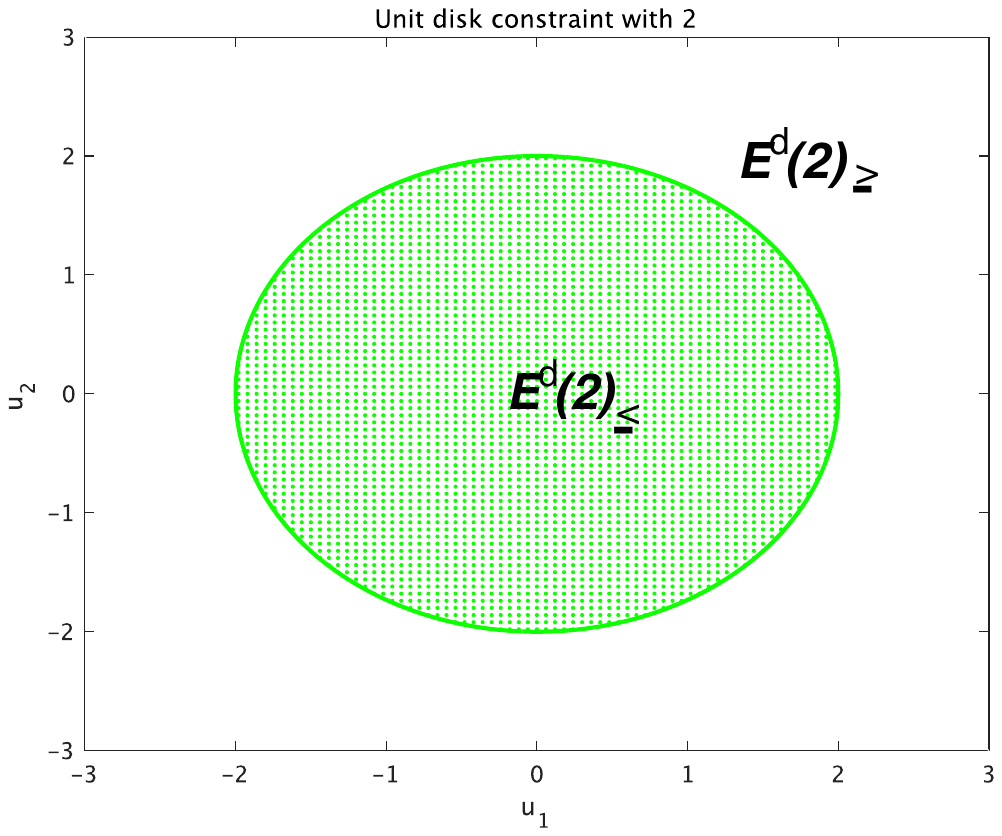}
\end{center}
\vspace{-35mm} 
\caption{The disk constraint $\disk(r)_\geq$: $r=1$ (left), $r=2$ (right).
} 
\end{figure}

\noindent
{\bf A  hyperbola constraint:}
\begin{eqnarray*}
& &\hyperbola(r) \equiv \begin{pmatrix} -1 & 0 & 0 \\ 0 & 1 & 0 \\ 0 & 0 & r^2 \end{pmatrix}, \
\hyperbola(r)_{\geq}  \ \mbox{or } \hyperbola(r)_{\leq} 
 = \left\{ \u \in \Real^2 : -u_1^2 + u_2^2 +r^2 \geq  \ \mbox{or } \leq 0 \right\},  
\end{eqnarray*}
where $r \geq 0$ denotes a parameter. See Figure 2.

\medskip

\begin{figure}[t!]   \vspace{-30mm} 
\begin{center}
\includegraphics[height=110mm]{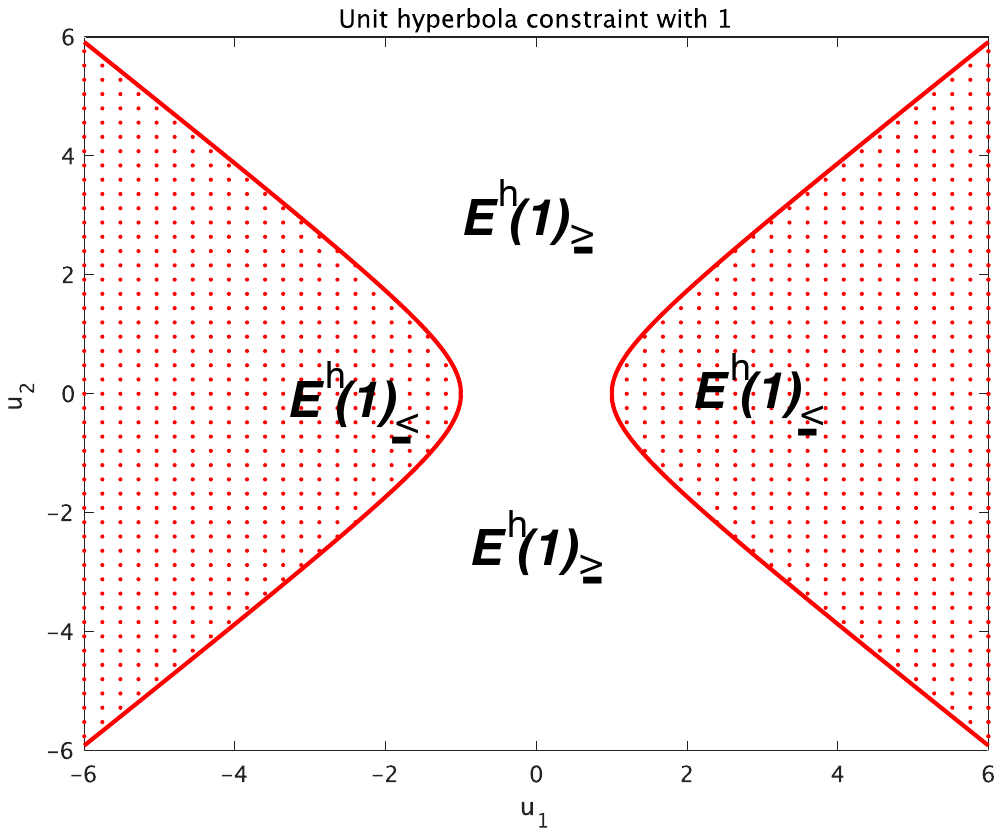}
\includegraphics[height=110mm]{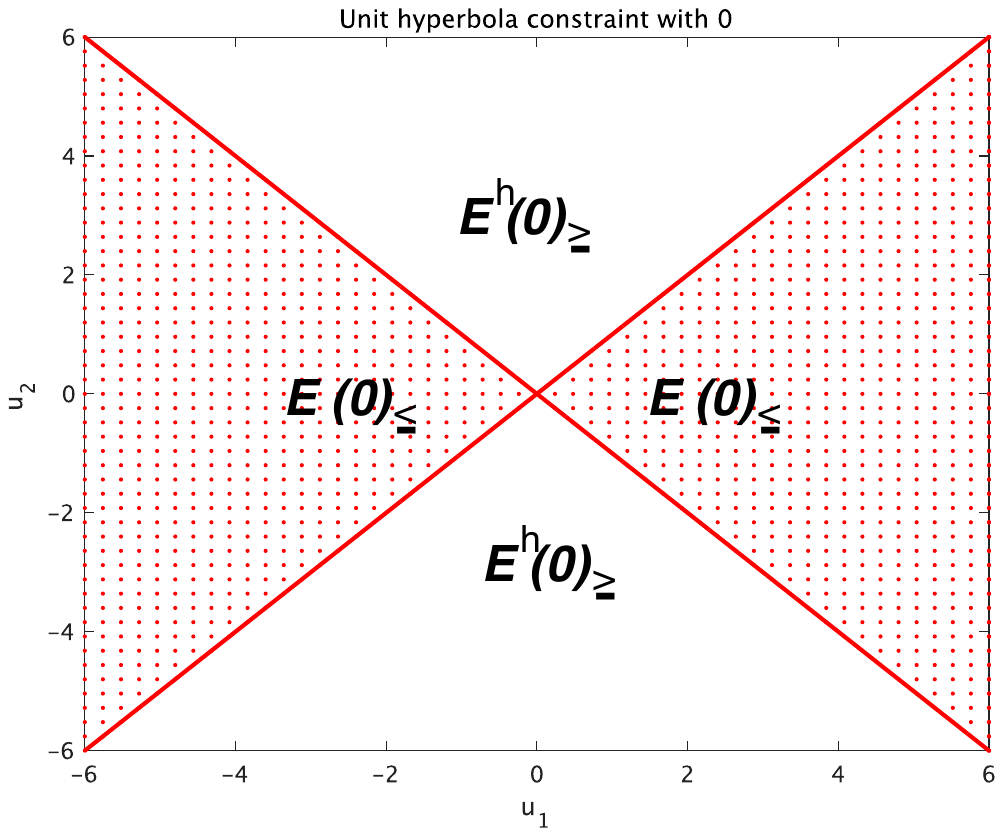}
\end{center}

\vspace{-35mm} 

\caption{
The hyperbola constraint $\hyperbola(r)_\geq$:  $r=1$ (left), $r=0$ (right).
}
\end{figure}  

\noindent
{\bf A parabola constraint: }
\begin{eqnarray*}
& &\parabola(r) \equiv \begin{pmatrix} 0 & 0 & -1/2 \\ 0 & 1 & 0 \\ -1/2 & 0 & r \end{pmatrix}, \
\parabola(r)_{\geq}  \ \mbox{or } \parabola(r)_{\leq} 
 = \left\{ \u \in \Real^2 : -u_1 + u_2^2 + r \geq   \ \mbox{or } \leq 0\right\},   
\end{eqnarray*}
where $r \geq 0$ denotes a parameter. See Figure 3.

\medskip

\begin{figure}[t!]   \vspace{-30mm} 
\begin{center}
\includegraphics[height=110mm]{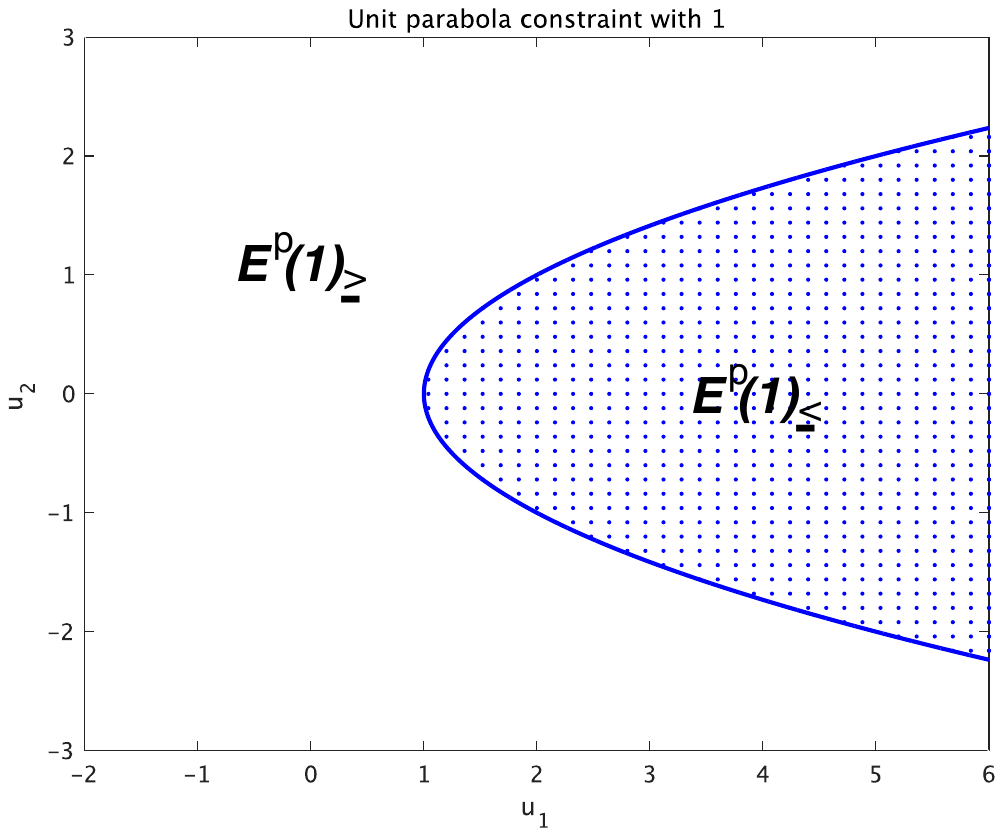}
\includegraphics[height=110mm]{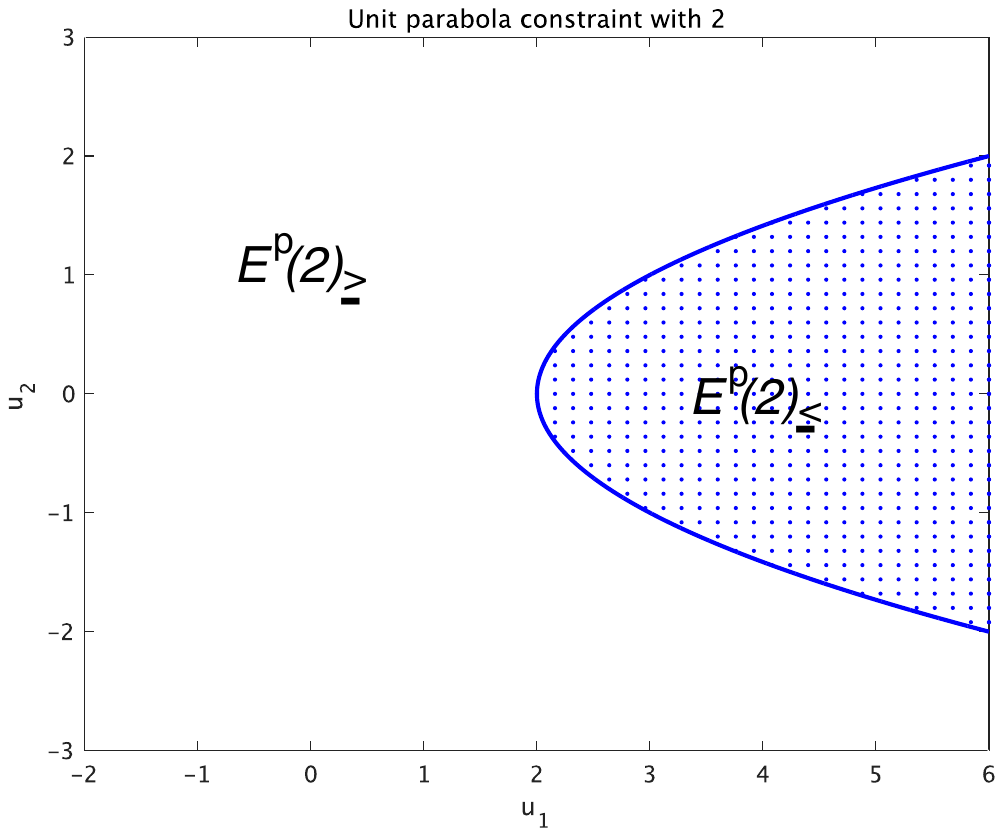}
\end{center}

\vspace{-35mm} 

\caption{
The parabola constraint $\parabola(r)_\geq$: $r=1$ (left), $r=2$ (right). 
}
\end{figure}

\noindent
{\bf A  linear constraint: }
\begin{eqnarray*}
& &\linear(r) \equiv \begin{pmatrix} 0 & 0 & 1/2 \\ 0 & 0 & 0 \\ 1/2 & 0 & -r \end{pmatrix}, \
\linear_{\geq}(r)  \ \mbox{or } \linear(r)_{\leq} 
 = \left\{ \u \in \Real^2 : u_1 - r \geq  \ \mbox{or } \leq 0 \right\},  
\end{eqnarray*}
where $r \in \Real$ denotes a parameter. See Figure 4. 

\medskip

\begin{figure}[t!]   \vspace{-30mm} 
\begin{center}
\includegraphics[height=110mm]{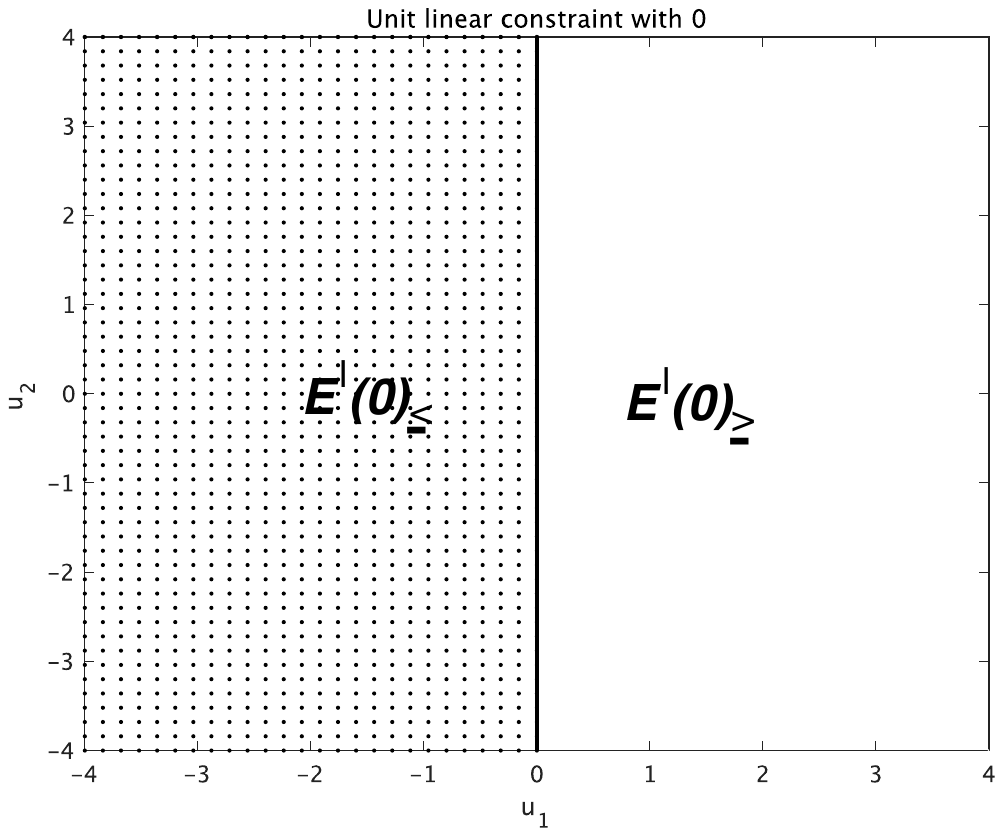}
\includegraphics[height=110mm]{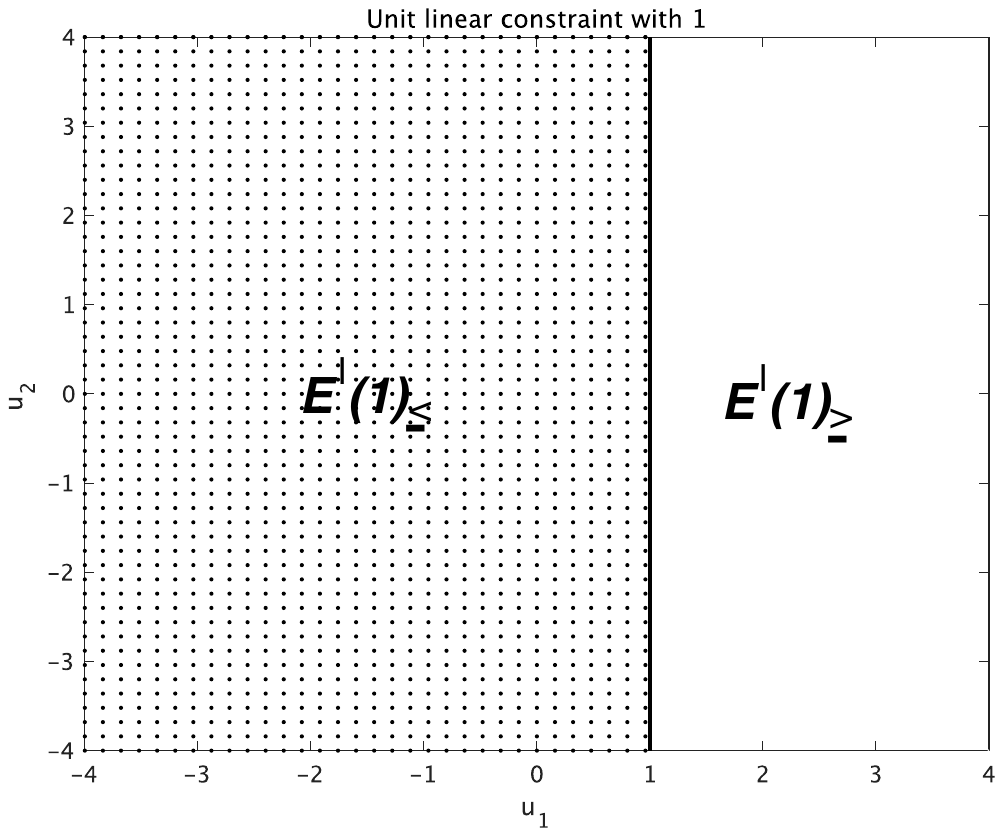}
\end{center}

\vspace{-35mm} 

\caption{The linear constraint  $\linear(r)_\geq$: $r=0$ (left), $r=1$ (right).}
\end{figure}  

\subsection{Scaling, rotation and parallel transformation}

If we apply a linear transformation defined by $ {\scriptsize \begin{pmatrix} \u \\ 1\end{pmatrix} } \rightarrow  
{\scriptsize \begin{pmatrix} \u \\ 1 \end{pmatrix}}=\M 
{\scriptsize \begin{pmatrix} \v \\ 1 \end{pmatrix}}$, where  ${\scriptsize \begin{pmatrix} \u \\ 1\end{pmatrix} } \in
\Real^3$ and  $\M$ denotes a $3 \times 3$ nonsingular  matrix, 
to the basic quadratic constraints described  in Section 2.1, 
we obtain general quadratic constraints in $\v$. 
Under this transformation, the quadratic function 
 ${\scriptsize\begin{pmatrix} \u \\ 1 \end{pmatrix}} ^T\B{\scriptsize \begin{pmatrix} \u \\ 1 \end{pmatrix}}$
in $\u \in \Real^2$ is transformed into the quadratic function  
 ${\scriptsize \begin{pmatrix} \v \\ 1 \end{pmatrix}} ^T\M^T\B\M
 {\scriptsize \begin{pmatrix} \v \\ 1 \end{pmatrix}}$. 
As  such a linear transformation,  we consider   a scaling, a rotation and a parallel 
 transformation, which are used in Section 2.3 for $3$ instances of $\BC$ satisfying 
Conditions (B)' and (D). 
 
\medskip

We describe matrices for scaling, rotation and parallel transformation.

\noindent
{\bf A $3\times 3$ scaling matrix: } 
$
\S(\ss) = \begin{pmatrix} 1/s_1 & 0 & 0 \\ 0 & 1/s_2 & 0 \\ 0 & 0 & 1 \end{pmatrix}
$, 
where $\0 < \ss =(s_1,s_2) \in\Real^2$ denotes a parameter vector. The vector $\u \in \Real^2$ 
is transformed to $\v = (s_1u_1,s_2u_2) \in \Real^2$. 

\medskip

\noindent
{\bf A $3\times 3$ rotation matrix:} 
$ 
\R(\theta) = \begin{pmatrix} \mbox{cos}\theta & \mbox{sin}\theta & 0 \\ -\mbox{sin}\theta & \mbox{cos}\theta & 0 \\ 0 & 0 & 1 \end{pmatrix}
$, 
where $\theta \in [-\pi,\pi]$ or $[0,2\pi]$ denotes a parameter. The vector  $\u \in \Real^2$ is transformed to
 $\v = (\mbox{cos}\theta u_1- \mbox{sin}\theta u_2,\mbox{sin}\theta u_1 
+ \mbox{cos}\theta u_2) \in \Real^2$. 

\medskip

\noindent
{\bf A $3\times 3$ matrix for parallel transformation:} 
$
\P(\p) = \begin{pmatrix} 1 & 0 & -p_1 \\ 0 & 1 & -p_2 \\ 0 & 0 & 1 \end{pmatrix}
$, 
\noindent
where $\p = (p_1,p_2) \in \Real^2$ denotes a parameter vector. 
The vector  $\u \in \Real^2$ is transformed to $\v = (u_1+p_1,u_2 + p_2) \in \Real^2$. 

\medskip

Figure 5 illustrates the application of scaling, rotation and  parallel transformation to 
the parabola constraint $\parabola(r)_\geq$ presented  in Section 2.1. 
Figure 6 illustrates the application of scaling, rotation and parallel transformation to 
the linear constraint $\linear(r)_\geq$ given in Section 2.1.
We note that if  the order of scaling, rotation and parallel transformation is changed, then
the resulting constraint may differ. 

\begin{figure}[t!]   \vspace{-30mm} 
\begin{center}
\includegraphics[height=110mm]{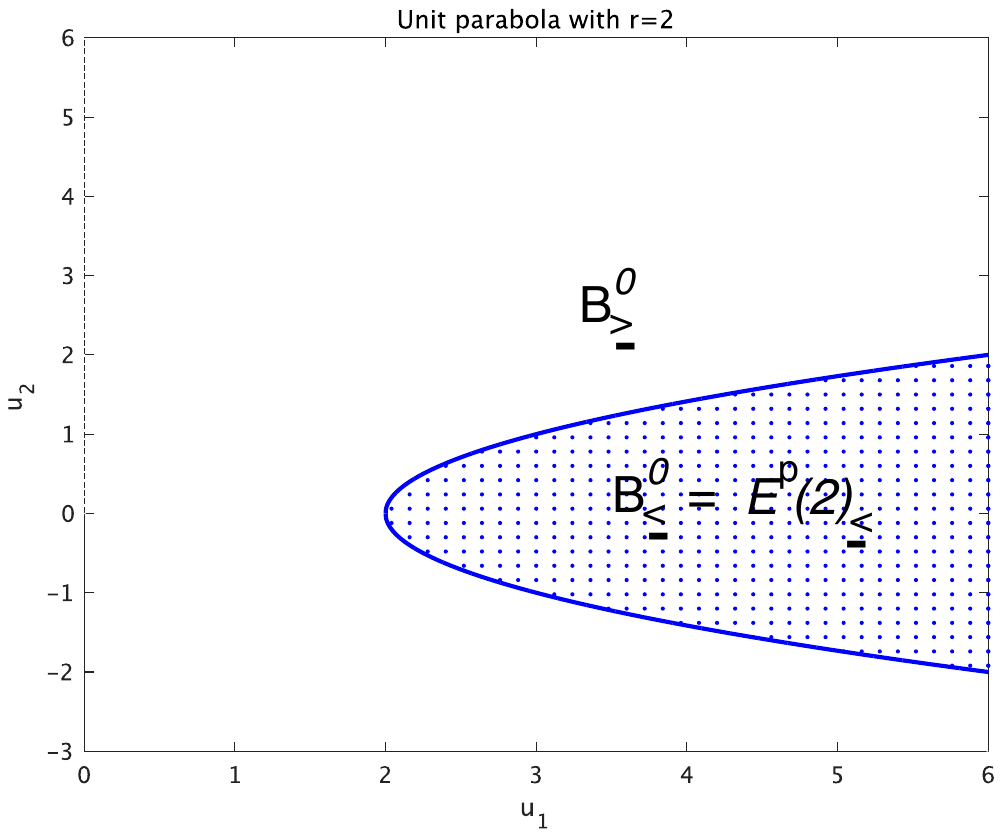}
\includegraphics[height=110mm]{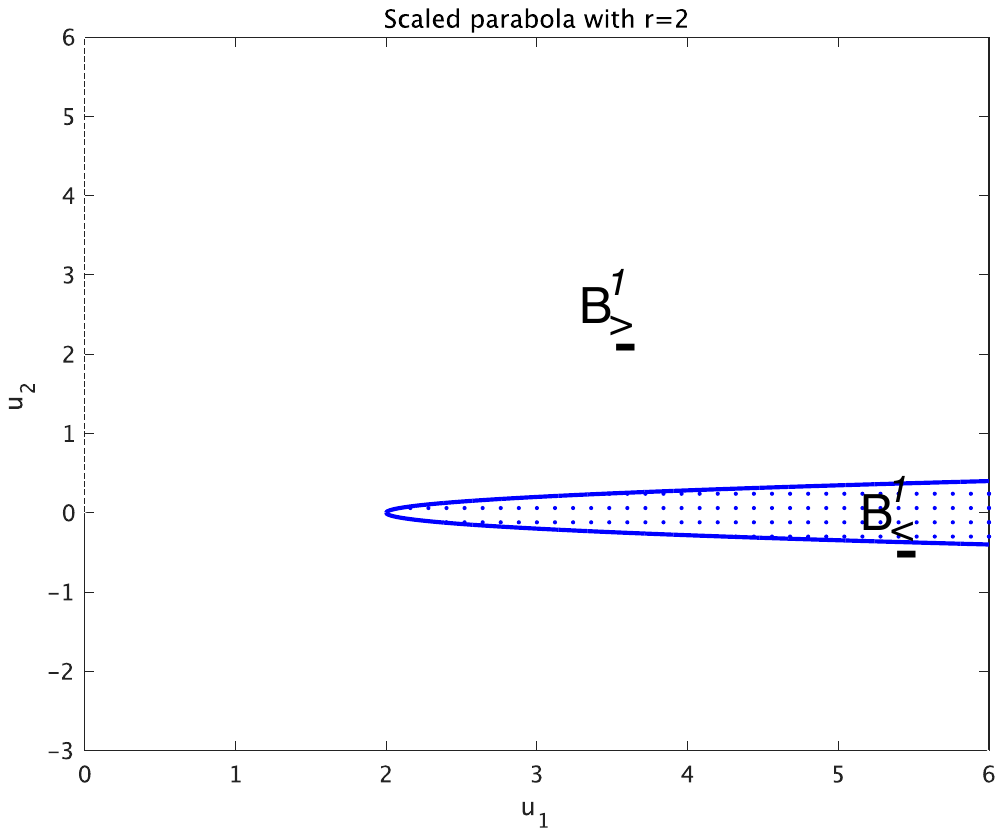}

\vspace{-55mm} 

\includegraphics[height=110mm]{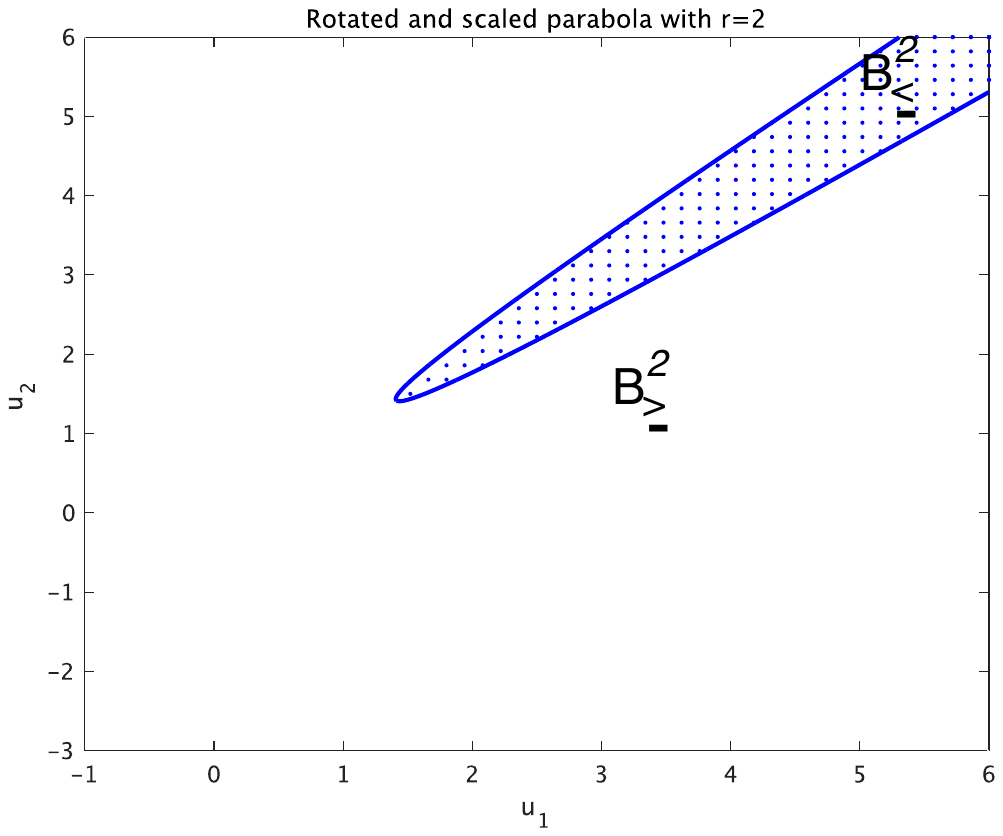}
\includegraphics[height=110mm]{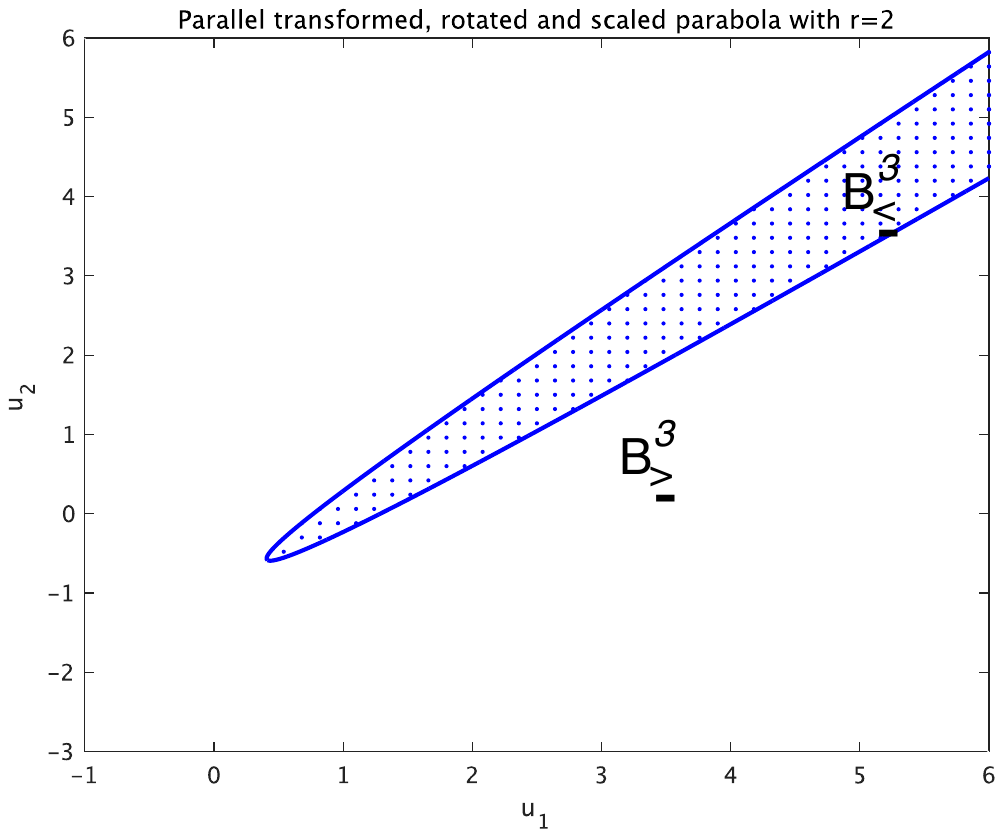}
\end{center}
\vspace{-35mm} 
\caption{The application of scaling, rotation and linear transformation to 
the parabola constraint $\parabola(r)_\geq$ described  in Section 2.1. Here 
$ 
 \B^0 = \mbox{$\parabola(2)$}, \  \B^1 = \S((1,0.2))^T \B^0 \S((1,0.2)), \ 
\B^2 = \R(\pi/4)^T \B^1 \R(\pi/4), \  \B^3 = \P((-1,-2))^T \B^2 \P((-1,-2)) 
$.  
}
\end{figure}  

\begin{figure}[t!]   
\begin{center}
\includegraphics[height=110mm]{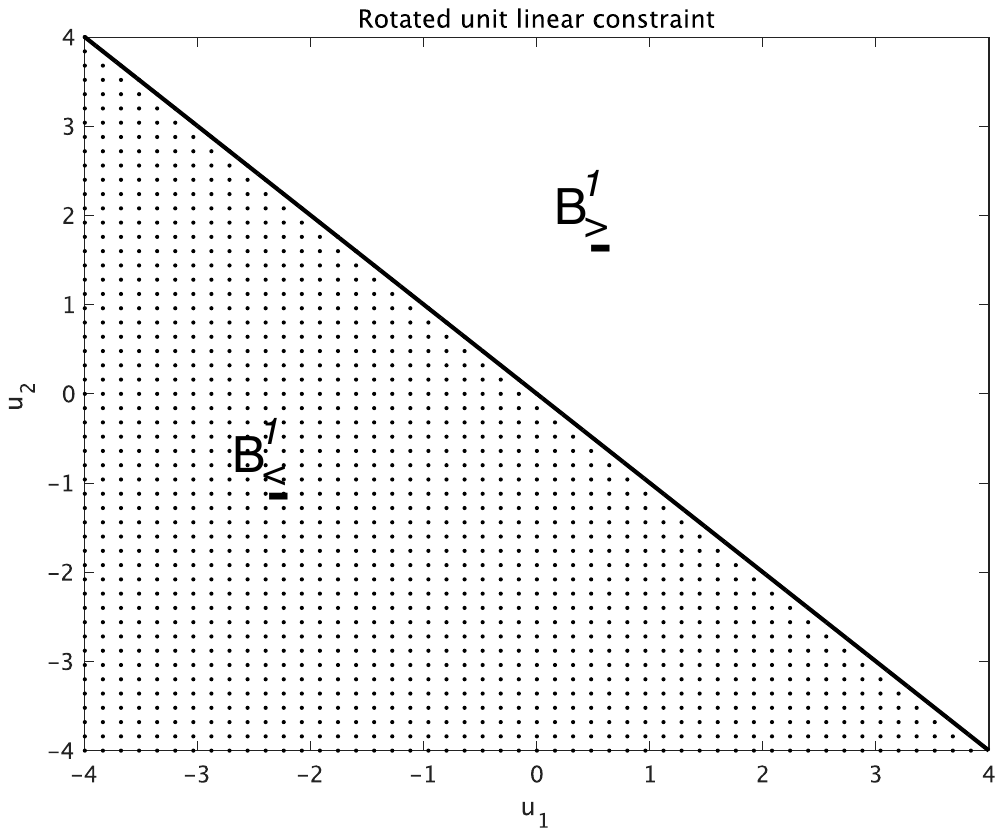}
\includegraphics[height=110mm]{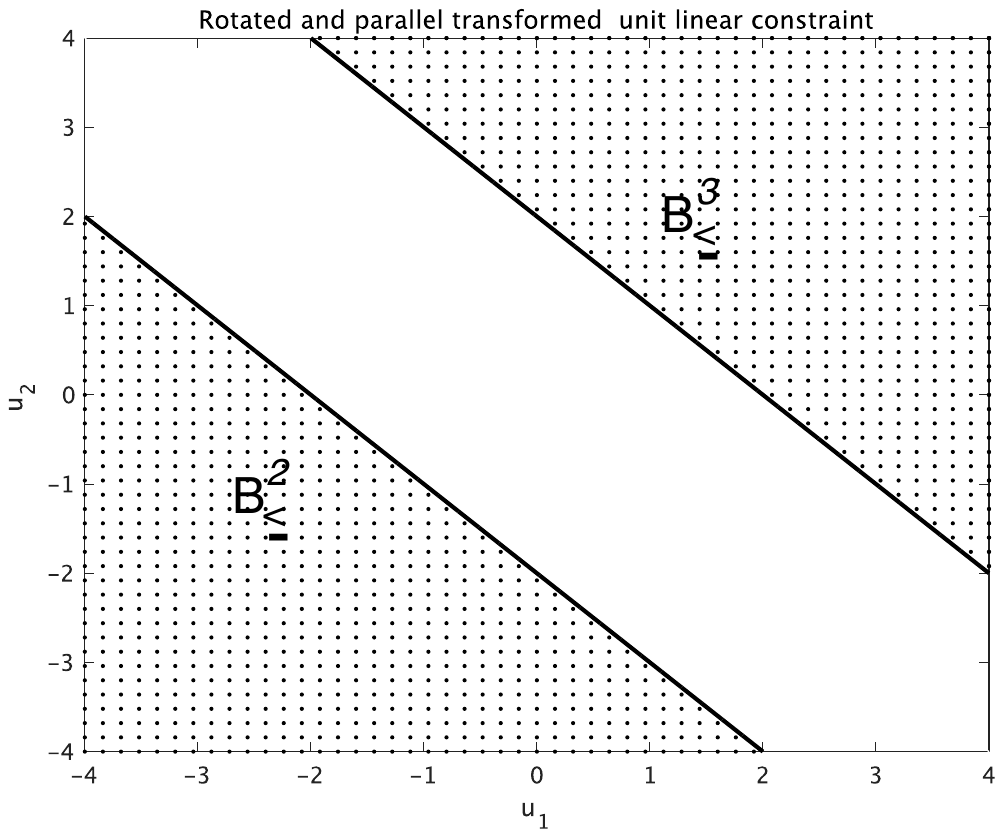}
\end{center}

\vspace{-35mm} 

\caption{
The application of scaling, rotation and  parallel transformation to 
the linear constraint $\linear(r)_\geq$
given in Section 2.1. 
  Here 
$\B^1 = \R(\pi/4)^T \mbox{$\linear(0)$}\R(\pi/4), \ 
\B^2 = \P((-1,-1))^T\B^1\P((-1,-1)), \ 
 \B^3 = \P((1,1))^T\R(5\pi/4)^T \mbox{$\linear(0)$}\R(5\pi/4)\P((1,1)). 
$ 
}
\end{figure}  


\subsection{
Instances of $\BC \subseteq \SymMat^3$ satisfying Conditions (B)' and (D)
}

Each $\B \in \BC \subseteq \SymMat^3$ in Instances 2.1, 2.2 and 2.3 is 
obtained by applying scaling,  rotation and/or 
 parallel transformation  to $\disk(r)$, $\hyperbola(r)$ or $\parabola(r)$ for some $r > 0$. 
Each $\B \in \BC \subseteq \SymMat^3$ in Instances 2.4 and 2.5 is described explicitly without 
relying on scaling, rotation or parallel transformation. 
Instance 2.6 leads to an extension of $\BC \subseteq \SymMat^3$ satisfying  Condition (D)  
to a general $\BC \subseteq \SymMat^n$ with $n \geq 3$ in Section 3. 
Throughout all instances of $\BC$, 
$\B_{\leq}$ is illustrated by the shaded regions,
 where the unshaded regions correspond to the feasible regions $\BC_\geq$ 
of QCQP \eqref{eq:QCQP3}. 
Instance 2.7 shows a case where $-\infty = \eta < \zeta$ occurs even if $\BC$  satisfies 
 Condition (B)'.

\examp \label{example:disk}
We define 
$\B^k \in \SymMat^3 \ (0 \leq k \leq 7)$ $(0 < r \leq 1/2)$ and  
$\BC \subseteq \SymMat^3$ as follows: 
\begin{eqnarray}
\left.
\begin{array}{l}
 \B^k = \R(k\pi/3)^T\P((1,0))^T\mbox{$\disk(r)$} \P((1,0)) \R(k\pi/3) \ (0 \leq k \leq 5), \\[3pt]
 \B^6 = \mbox{$\disk(r)$}, \  \B^7 = -\disk(3/2),\
 \BC = \{\B^k : 0 \leq k \leq 7 \}. 
\end{array}
\right\} \label{eq:disk}
\end{eqnarray}
Figure 7 illustrates 
$\B^k_\leq$ $(0 \leq k \leq 7)$, where the unshaded region 
corresponds to $\BC_\geq$. We took $r=1/2$ in the left figure, and $r = 1/3$ in the 
right figure. We see from Figure 7 that $\BC$ satisfies  Condition (B)'. 
Let $\alpha_{B^k} = 1$ $(0 \leq k \leq 6)$ and $\alpha_{B^7} = 1/3$. We can verify that 
$\alpha_{B^j}\B^j + \alpha_{B^k}\B^k \in \SymMat^3_+$ $(0 \leq j < k \leq 7)$. Therefore, 
$\BC$  also satisfies Condition (D).
\eexamp

\begin{figure}[t!]   \vspace{-30mm} 
\begin{center}
\includegraphics[height=110mm]{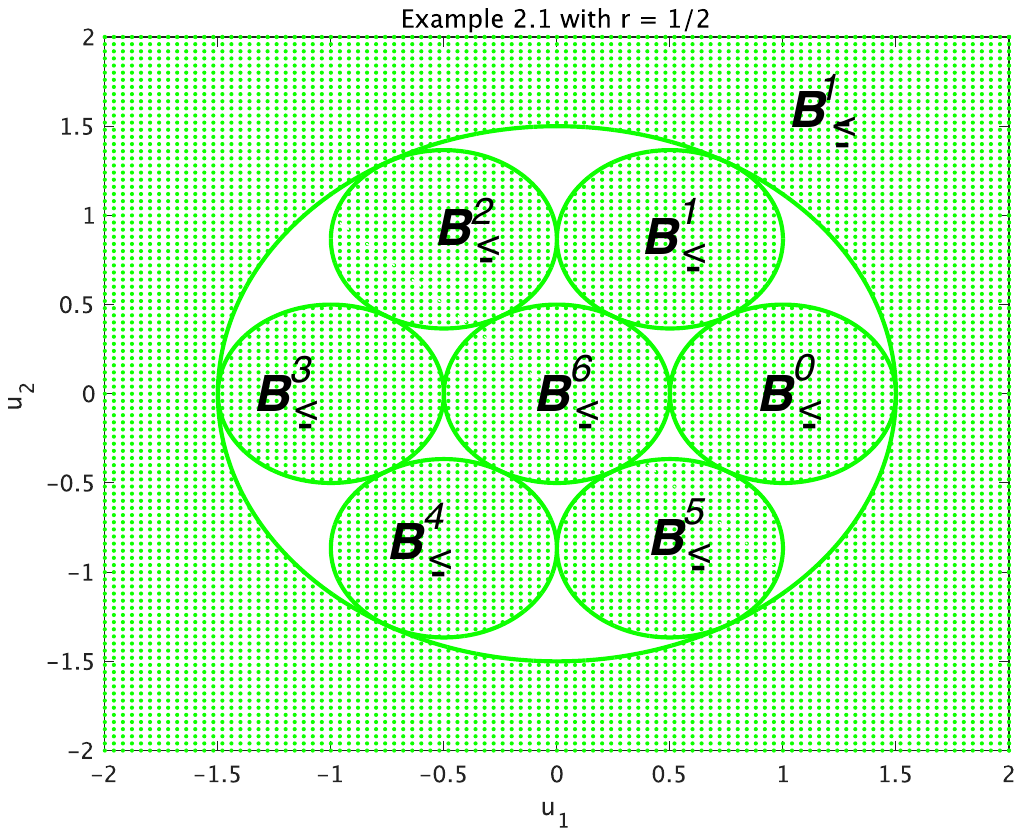}
\includegraphics[height=110mm]{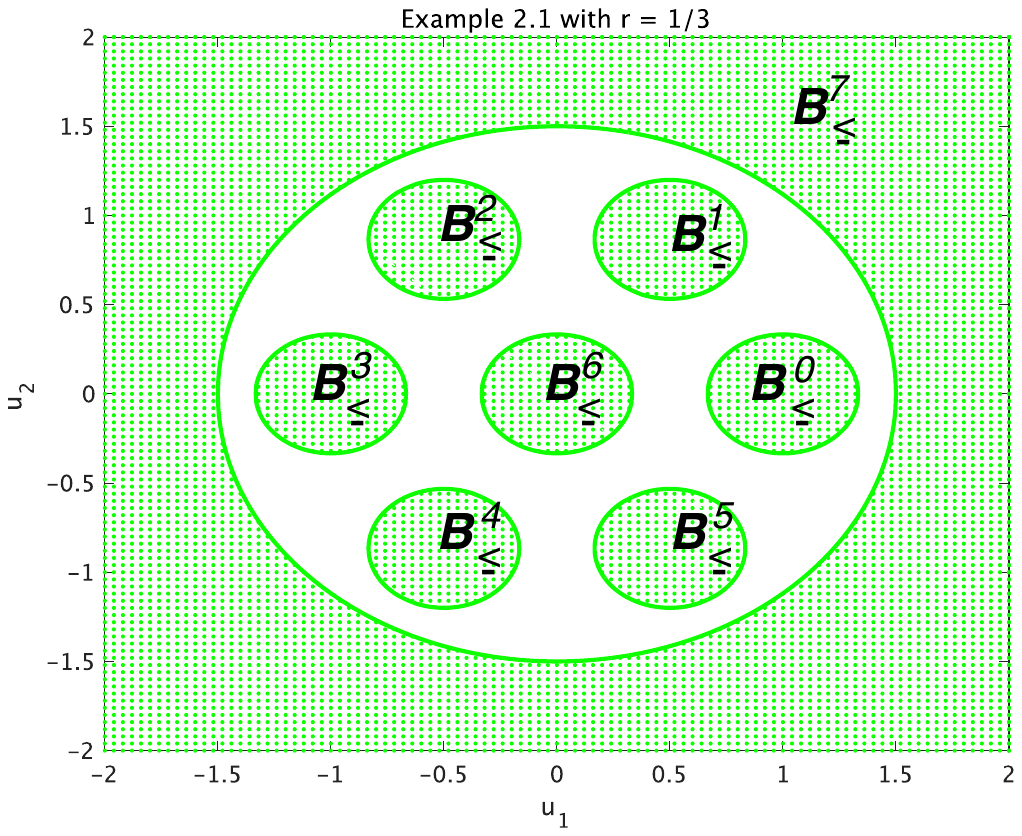}
\end{center}
\vspace{-30mm} 
\caption{
Illustration of Instance 2.1. The value of $r$ is $1/2$ for the left figure and  $1/3$ for the right figure.
}
\end{figure}  

\examp 
\label{example:hyperbola}
Let $2 \leq m \in \Integer$ where $\Integer$  denotes the set of integers. We also let $r > 0$ and $\p \in \Real^2$. 
Define
\begin{eqnarray*}
& &  \B^k = \P(\p)^T
\R(k\pi/m)^T\S((1,\mbox{tan}(\pi/(2m))))^T{\hyperbola(r)}\\
& & \hspace{29mm} \S((1,\mbox{tan}(\pi/(2m))))\R(k\pi/m)\P(\p) \ (0 \leq k \leq m-1), \\
& & \B^m = \P(\p)^T\mbox{$\disk(r)$}\P(\p), \
\BC = \left\{  \B^k \ (0 \leq k \leq m) \right\} . 
\end{eqnarray*}
See Figure 8. In the left figure, we took $m=2$, $r=1$ and $\p =(1,1)$, where the  hyperbola constraint 
\begin{eqnarray*}
\B_\leq & = & (\S((1,\mbox{tan}(\pi/(2\times2))))^T{\hyperbola(1)}
\S((1,\mbox{tan}(\pi/(2\times2)))))_\leq \\
 &  = & 
\left\{ \u \in \Real^2 : -u_1^2 + u_2^2 + 1^2 \leq 0 \right\}= \hyperbola(1)_\leq
\end{eqnarray*} 
is shifted to the upper-right direction with the origin $(0,0) \rightarrow (1,1)$ to 
create  $\B^0_\leq$. 
After rotating $\B_\leq$ by $\pi/2$, and then we move the resulting hyperbola toward the 
$(1,1)$-direction to obtain $\B^1_\leq$. 
We can make a similar observation on the right figure. From Figure 8, we see that 
$\BC$ satisfies Condition (B)'. 
We  can also verify that 
$\BC$ also satisfies  Condition (D) with $\alpha_{B^k} = 1$ $(0 \leq k \leq m)$. 
\eexamp

\begin{figure}[t!]   \vspace{-30mm} 
\begin{center}
\includegraphics[height=110mm]{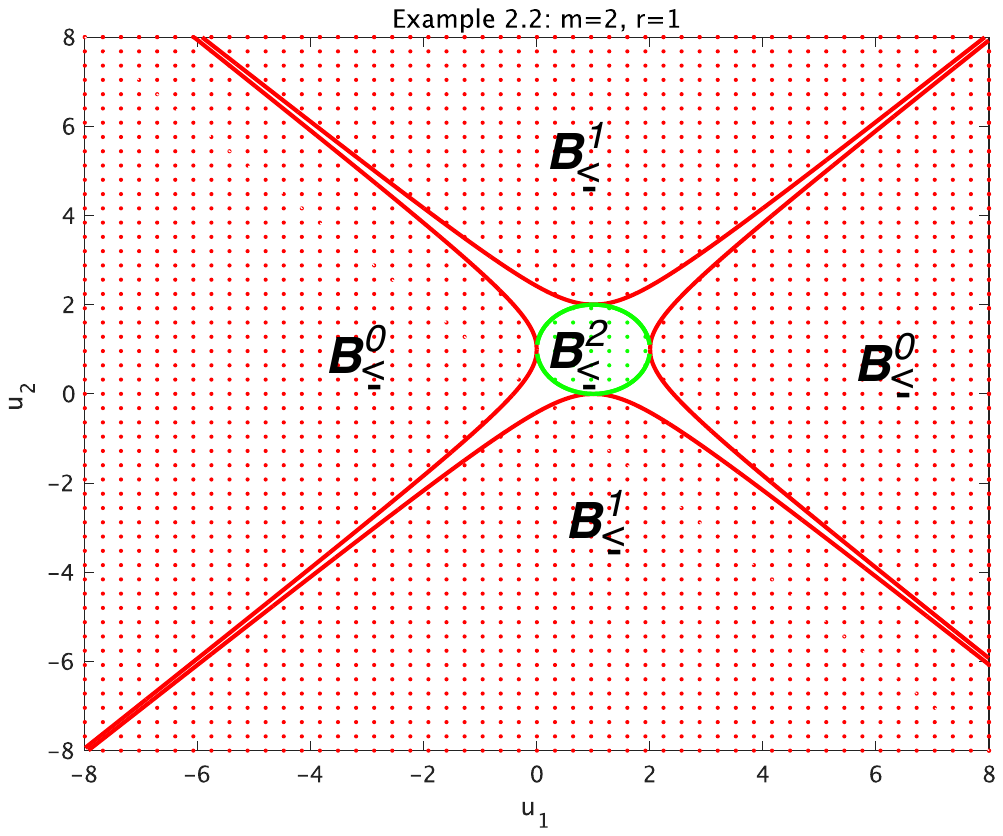}
\includegraphics[height=110mm]{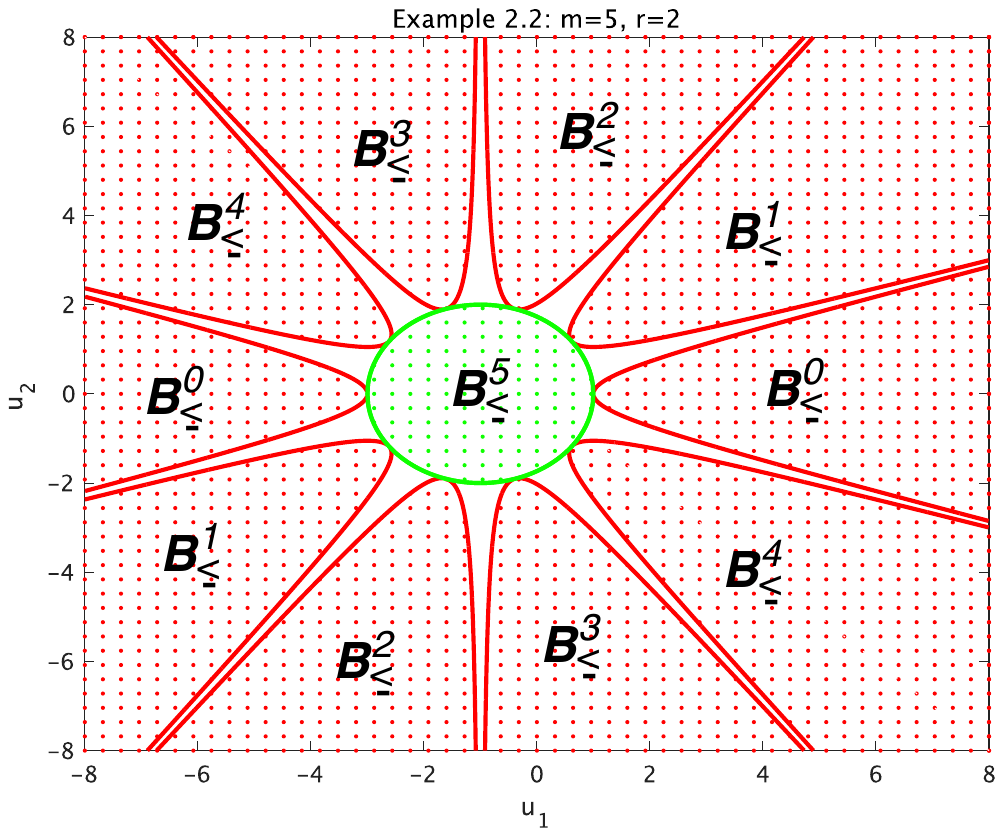}
\end{center}
\vspace{-30mm} 
\caption{Illustration of Instance 2.2. The values for the parameters are  $m=2,\ r=1$ and $\p = (1,1)$ for the left figure, and 
$m = 5, \ r = 2$ and $\p = (-1,0)$ for the right figure.
}
\end{figure}  

\examp \label{example:parabola}
For every $3 \leq m \in \Integer$ and $r > 0$,  
define 
\begin{eqnarray}
\left. 
\begin{array}{l}
 \C^k = 
\R(2k\pi/m)^T\S((1,2(\mbox{tan}(\pi/m))\sqrt{r}))^T{\parabola(r)}\\
\mbox{ \ }  \hspace{43mm} \S((1,2(\mbox{tan}(\pi/m))\sqrt{r})) \R(2k\pi/m) \ (0 \leq k \leq m-1),  \\[3pt]
\C^m =  \mbox{$\disk(r)$}, \ \CC = \left\{  \C^k \ (0 \leq k \leq m)\right\}.
\end{array}
\right\} \label{eq:parabola} 
\end{eqnarray}
See Figure 9.  In the left figure, $\C^0_\leq$ is expressed as a scaled parabola constraint such that 
\begin{eqnarray*}
\C^0_\leq & = & (\S((1,2(\mbox{tan}(\pi/3))\sqrt{1}))^T\mbox{$\parabola(1)$}
\S((1,2(\mbox{tan}(\pi/3))\sqrt{1})))_\leq \\
 &  = & 
\left\{ \u \in \Real^2 : -u_1 + ((2\mbox{tan}(\pi/3))u_2)^2 + 1 \leq 0 \right\}. 
\end{eqnarray*}
$\C^1_\leq$ and $\C^2_\leq$ 
are obtained by 
rotating  $\C^0_\leq$ by $\pi/3$ and $2\pi/3$, respectively. 
We can make a similar observation on the right figure. 
We see from Figure 9 that $\CC$ satisfies  Condition (B)'. 
Let $\alpha_{B^k} = 1$ $(0 \leq k \leq m-1)$ and $\alpha_{B^m} = 1/(2r)$. Then 
we  can verify that 
$\alpha_{C^j} \C^j + \alpha_{C^k} \C^k \in \SymMat^n_+$
$(0 \leq j < k \leq m)$ holds. 
Therefore,  $\CC$ satisfies  Condition (D). 
\eexamp

\begin{figure}[t!]   \vspace{-30mm} 
\begin{center}
\includegraphics[height=110mm]{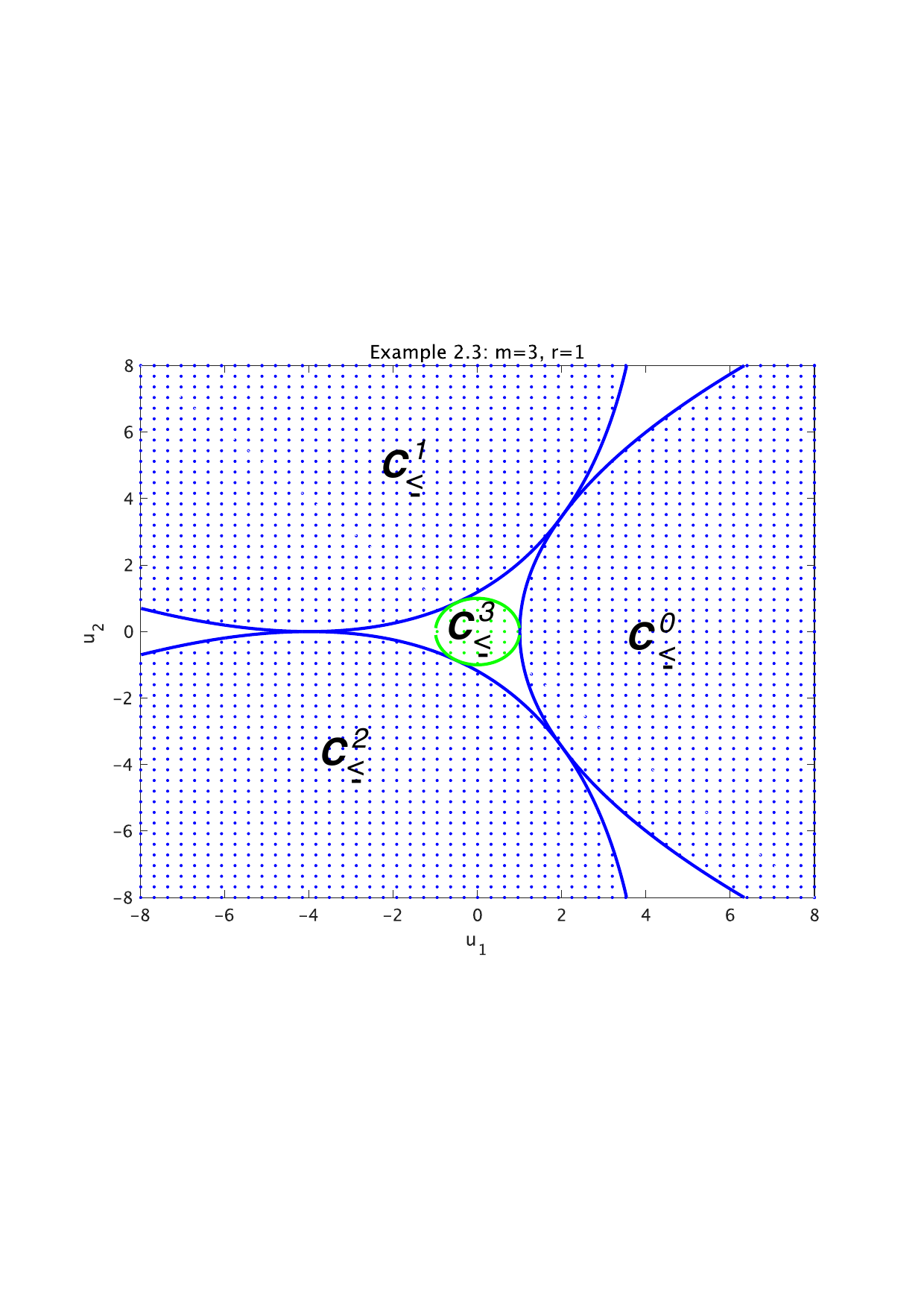}
\includegraphics[height=110mm]{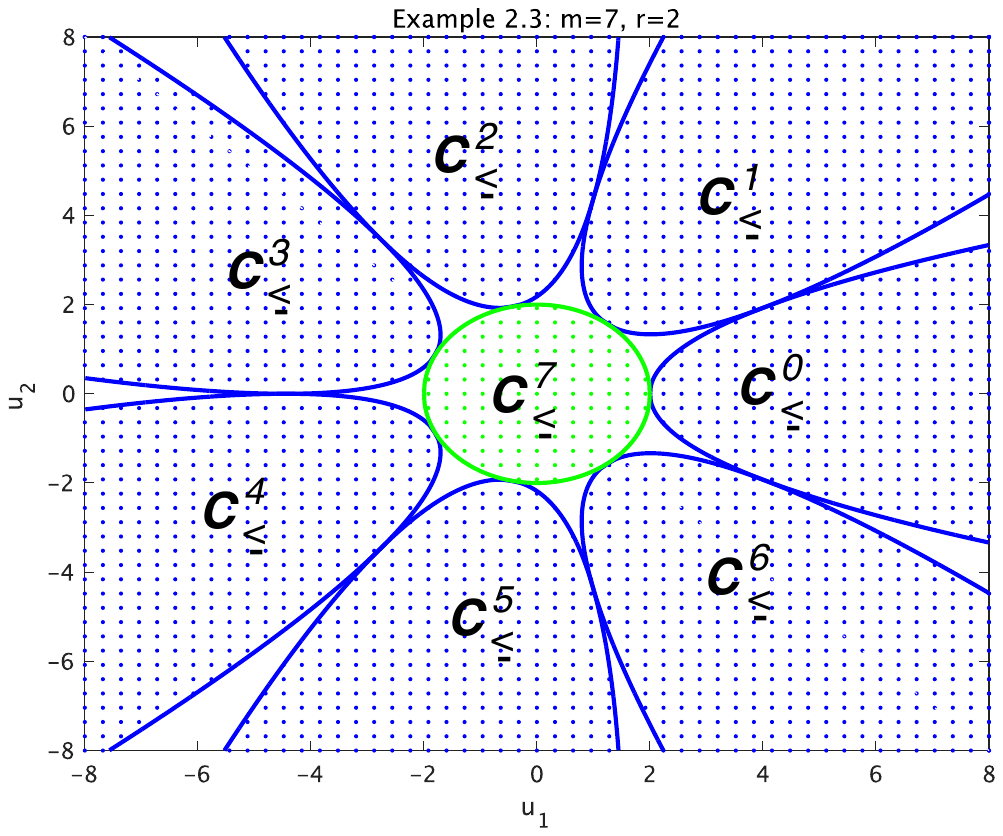}
\end{center}
\vspace{-30mm} 
\caption{Illustration of Instance 2.3. We took $(m,r)=(3,1)$ for the left figure, and $(m,r)=(7,2)$ for the right figure}. 
\end{figure} 

\examp \label{example:hyperbolaInf}
In this and next instances, we construct $\BC$ without relying on scaling, rotation and 
parallel transformation. 
For every $a \in \Integer$ and $r \geq 0$, define 
\begin{eqnarray*}
& & \B(a,r) = \begin{pmatrix} a^2 - 1/4& -a&0\\ -a &1&0\\ 0&0&r^2 \end{pmatrix}. 
\end{eqnarray*}
We then see that 
\begin{eqnarray*}
 \B(a,r)_\geq \ \mbox{or } \B(a,r)_\leq & = & 
\left\{ \u \in \Real^2 : 
(u_2 -au_1)^2 - u_1^2/4 + r^2 \geq \ \mbox{or } \leq 0 \right\}. 
\end{eqnarray*}
It is easy to observe that, for every distinct $a_1, \ a_2 \in \Integer$ and every $r_1, \ r_2 \in \Real_+$,  
\begin{eqnarray*}
\B(a_1,r_1) + \B(a_2,r_2) & = & 
\begin{pmatrix} a_1^2 + a_2^2 - 1/2&-a_1-a_2& 0\\ -a_1-a_2 &2 &0\\ 0&0&r_1^2 +r_2^2 \end{pmatrix}
\in \SymMat^3_+ 
\end{eqnarray*}
holds. 
Hence, if we take a finite subset ${G^{\rm h}}$ of $\Integer \times \Real_+$ satisfying 
\begin{eqnarray*}
a_1 \not= a_2 \ 
\mbox{for every pair of distinct $(a_1,r_1) \in {G^{\rm h}}$ and  $(a_2,r_2) \in {G^{\rm h}}$},  
\end{eqnarray*}
then $\BC({G^{\rm h}}) = \left\{ \B(a,r) : (a,r) \in {G^{\rm h}} \right\}$ satisfies  Condition (D). 
By Theorem \ref{theorem:main} (iv), $\BC({G^{\rm h}})$ also satisfies  Condition (B)', 
as can be verified by Figure 10.
\eexamp

\begin{figure}[t!]   \vspace{-30mm} 
\begin{center}
\includegraphics[height=110mm]{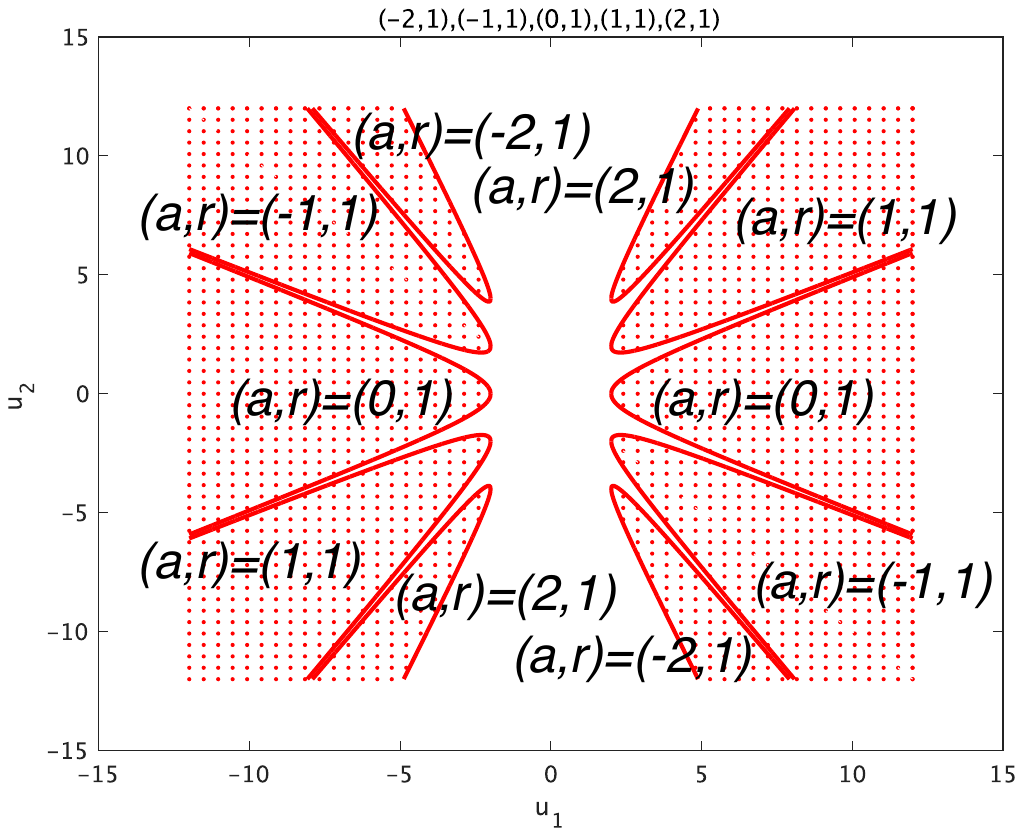}
\includegraphics[height=110mm]{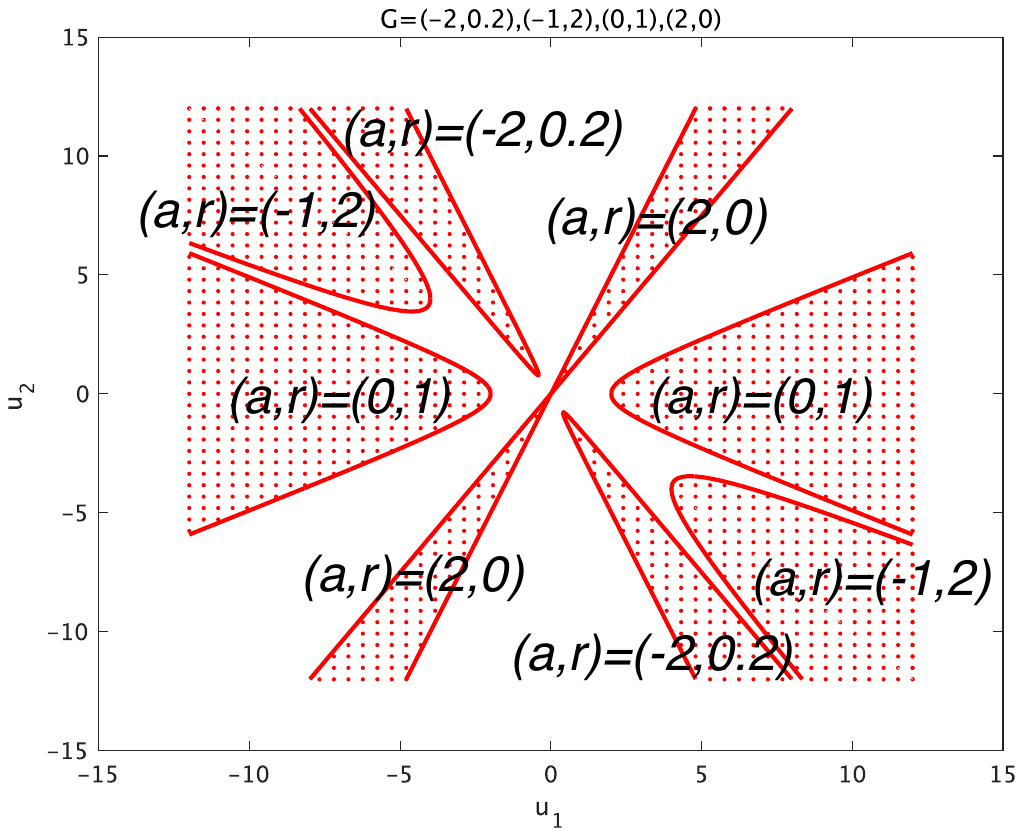}
\end{center}
\vspace{-30mm}
\caption{Illustration of  Instance 2.4: We took ${G^{\rm h}}=\{(2,1),(1,1),(0,1),(-1,1),(-2,1)\}$ for the left figure, 
and 
$
{G^{\rm h}} =\{(2,0),(0,1),(-1,2),(-2,0.2)\}$ for the right figure.
}
\end{figure}  

\examp \label{example:parabolaInf}
For every $a \in \Integer$ and $r \geq 0$, define 
\begin{eqnarray}
& & \B(a,r) = \begin{pmatrix} a^2 & -a & -1/2 \\ -a&1&0\\ -1/2 &0&r \end{pmatrix}.
\label{eq:parabolaInf}
\end{eqnarray}
Then it follows that 
\begin{eqnarray*}
 \B(a,r)_\geq \ \mbox{or } \B(a,r)_\leq & = & 
\left\{ \u \in \Real^2 : 
(au_1 - u_2)^2 - u_1 + r \geq 0 \ \mbox{or } \leq 0 \right\}.  
\end{eqnarray*} 
We can easily verify that, for every distinct $a_1, \ a_2 \in \Integer$ and every $r_1, \ r_2 \in [1,\infty)$, 
\begin{eqnarray*}
\B(a_1,r_1) + \B(a_2,r_2) & = & 
\begin{pmatrix} a_1^2 + a_2^2& -a_1-a_2& -1\\ -a_1-a_2 &2 &0\\ -1&0&r_1 +r_2 \end{pmatrix}
\in \SymMat^3_+
\end{eqnarray*}
holds. 
Therefore, for every finite subset ${G^{\rm p}}$ of  $\Integer \times [1,\infty)$ satisfying 
\begin{eqnarray*}
a_1 \not= a_2 \ 
\mbox{for every pair of distinct $(a_1,r_1) \in {G^{\rm p}}$ and $(a_2,r_2) \in {G^{\rm p}}$}, 
\end{eqnarray*}
$\BC({G^{\rm p}}) = \left\{ \B(a,r) : (a,r) \in {G^{\rm p}} \right\}$ satisfies  Condition (D). 
By Theorem \ref{theorem:main} (iv), $\BC({G^{\rm p}})$ satisfies  Condition (B)'. See Figure 11. 
\eexamp

\begin{figure}[t!]   \vspace{-30mm} 
\begin{center}
\includegraphics[height=110mm]{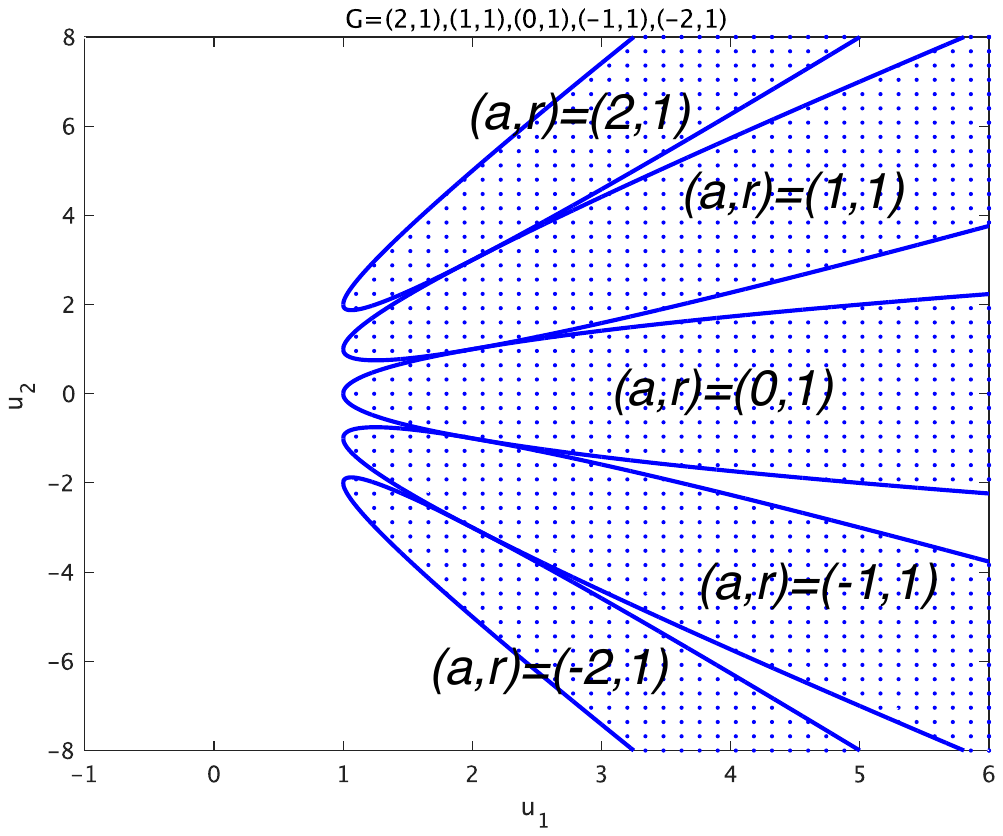}
\includegraphics[height=110mm]{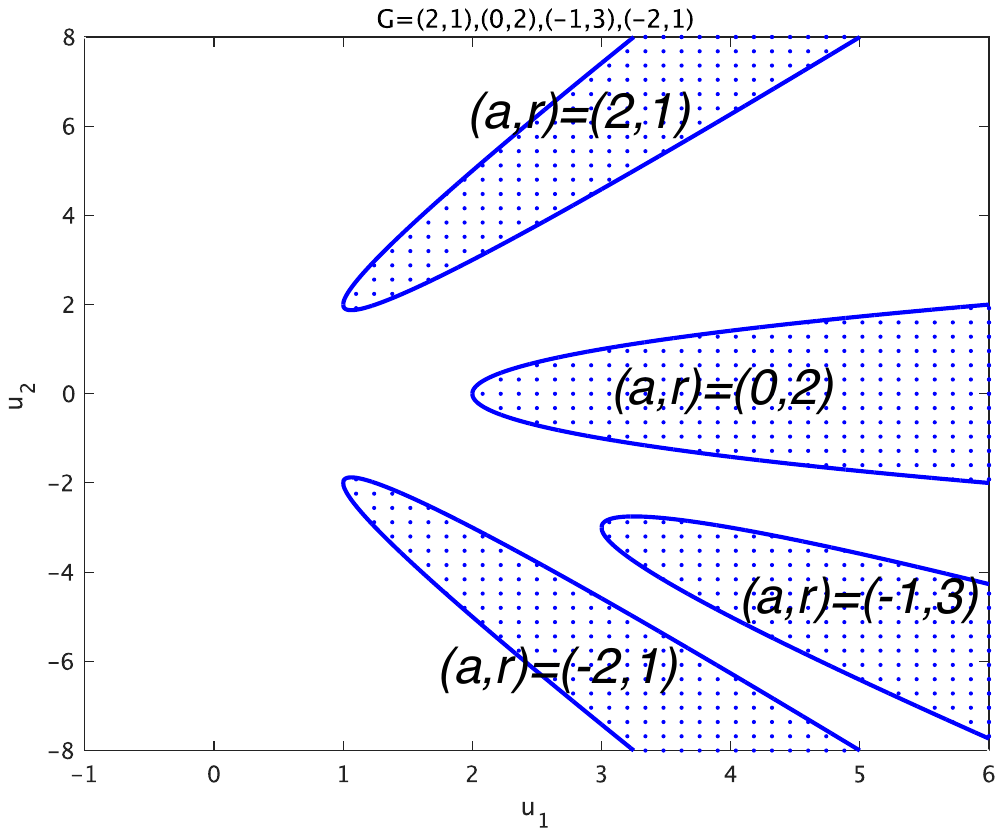}
\end{center}
\vspace{-30mm} 
\caption{Illustration of Instance~\ref{example:parabolaInf}.  We took ${G^{\rm p}}=\{(2,1),(1,1),(0,1),(-1,1),(-2,1)\}$ for the left figure, 
and 
${G^{\rm p}}=\{(2,1),(0,2),(-1,3),(-2,1)\}$ for the right figure.
}
\end{figure}  

\examp \label{example:combination}
In this instance, we utilize $\BC$  in \eqref{eq:disk} (Instance 2.1) and 
$\CC$  in \eqref{eq:parabola} (Instance 2.3). 
Let $r=1/2$ as in  the left figure of Figure 7, and $(m,r') = (7,2)$ as in the right figure of Figure 9. 
Then each of $\BC$ and $\CC$ consists of $8$ matrices in $ \SymMat^3$.
As we have stated there, they both satisfy  Condition (D); 
\begin{eqnarray*}
& & \alpha_{B^j}\B^j + \alpha_{B^k}\B^k \in \SymMat^3_+,\ \alpha_{B^j} > 0, \ \alpha_{B^k} >0 \ (0 \leq j < k \leq 7), \\
& & \alpha_{C^j}\C^j + \alpha_{C^k}\C^k \in \SymMat^3_+,\ \alpha_{C^j} > 0, \ \alpha_{C^k} >0 \ (0 \leq j < k \leq 7). 
\end{eqnarray*}
Letting $\lambda > 0$, $\mu > 0$ and $\lambda + \mu = 1$, we consider the convex 
combination of $\alpha_{B^k}\B^k$ and $\alpha_{C^k}\C^k$ such that 
\begin{eqnarray*}
\A^k = \lambda\alpha_{B^k}\B^k + \mu \alpha_{C^k}\C^k\ (0 \leq k \leq 7). 
\end{eqnarray*}
Then 
\begin{eqnarray*}
\A^j + \A^k & = & \lambda\alpha_{B^j}\B^j + \mu \alpha_{C^j}\C^j
+ \lambda\alpha_{B^k}\B^k + \mu \alpha_{C^k}\C^k\\ 
& = &  \lambda(\alpha_{B^j}\B^j + \alpha_{B^k}\B^k ) 
+ \mu (\alpha_{C^j}\C^j + \alpha_{C^k}\C^k) \in \SymMat^3_+ \ (0\leq j < k \leq 7)
\end{eqnarray*}
holds. 
Therefore $\AC = \{\A^k : 0 \leq k \leq 7\}$ satisfies  Condition (D). See Figure 12. 
If we renumber $\C^k$ contained in $\CC$ 
before taking the convex combination, a different $\AC$ is obtained. 
\eexamp

\begin{figure}[t!]   \vspace{-30mm} 
\begin{center}
\includegraphics[height=110mm]{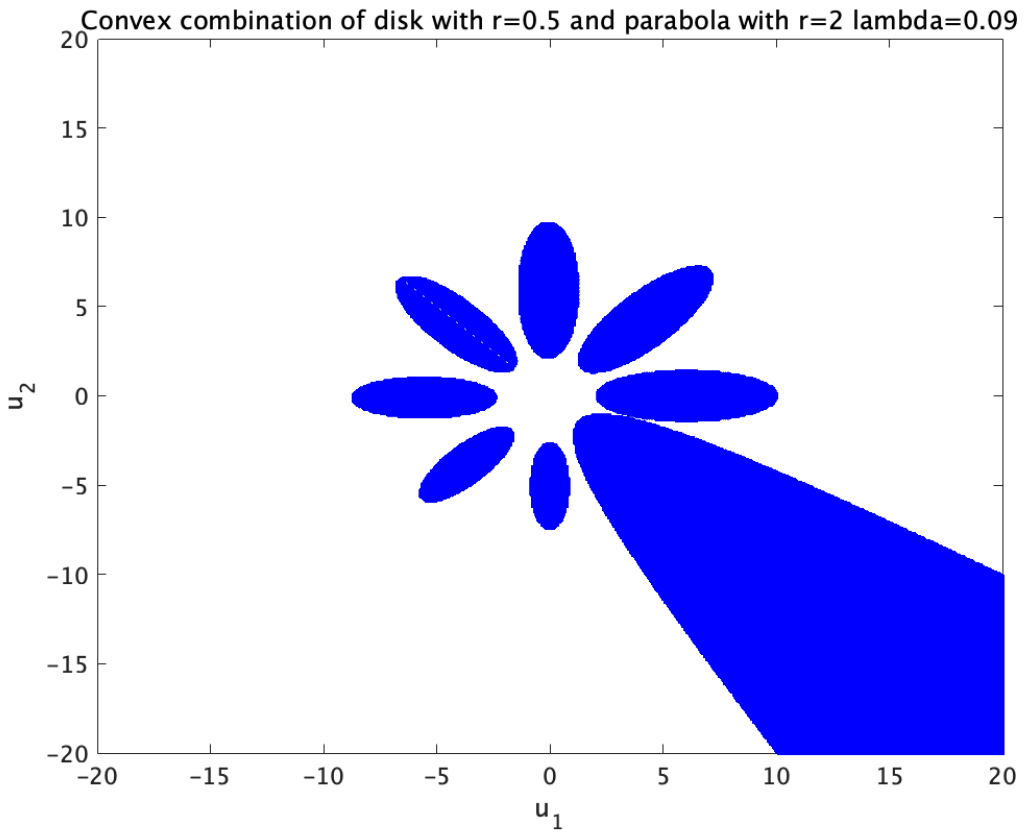}
\includegraphics[height=110mm]{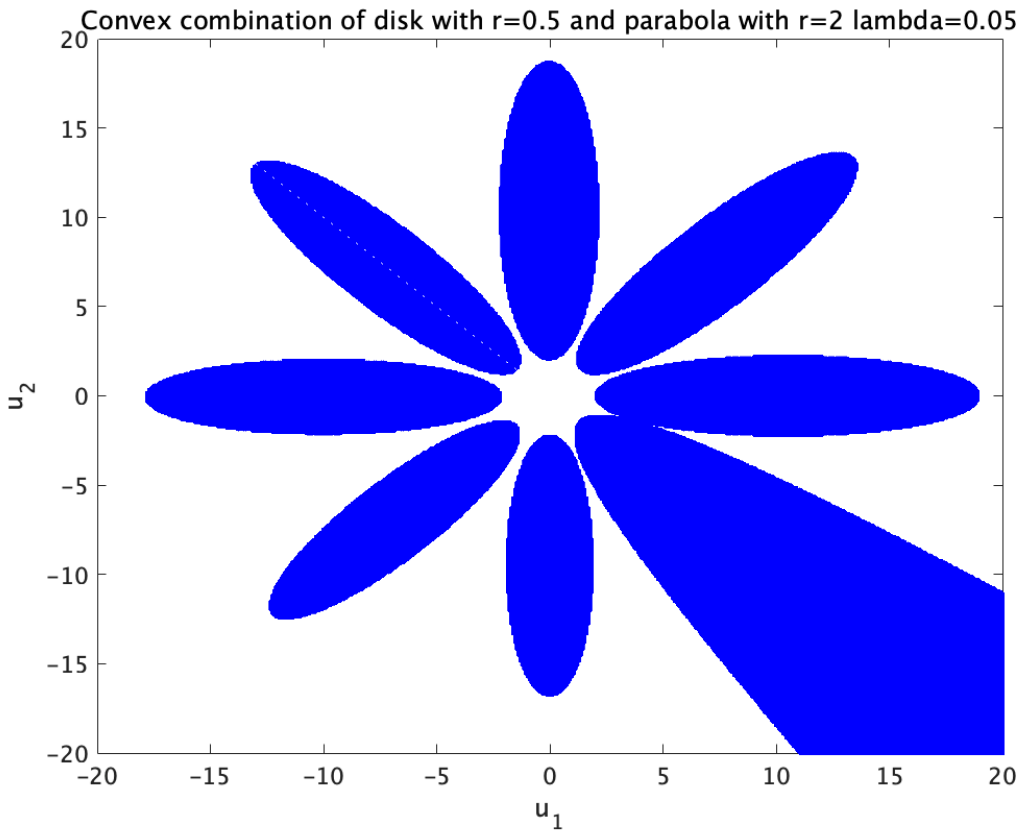}
\end{center}
\vspace{-30mm} 
\caption{Illustration of Instance 2.6: We took $(\lambda,\mu) = (0.09,0.91)$ for the left figure, and 
$(\lambda,\mu) = (0.05,0.95)$ for the right figure. 
}
\end{figure}  

\examp 
\label{example:counter}
Define 
\begin{eqnarray*}
& &
\B^2 = \begin{pmatrix} 0 & 0 & 1/2 \\  0 & 0 & 1/2 \\ 1/2 & 1/2 & 2 \end{pmatrix} \in \SymMat^3, \
 \B^3 = \begin{pmatrix} 0 & 0 & -1/2 \\  0 & 0 & -1/2 \\ -1/2 & -1/2 & 2 \end{pmatrix},  
\end{eqnarray*}
and let $\BC = \{\B^2, \ \B^3 \}$. Then it follows that 
\begin{eqnarray}
& & \B^2_\leq = \{ \u \in \Real^2 :  u_1 + u_2 + 2 \leq 0 \}, \ 
\B^3_\leq = \{ \u \in \Real^2 :  -u_1 - u_2 + 2 \leq 0 \}, \  \B^2_\leq \cap \B^3_\leq 
= \emptyset, \nonumber \\ 
& & \BC_\geq = \{ \u \in \Real^2 :  -2 \leq u_1+u_2 \leq 2 \}. \label{eq:linearConst}
\end{eqnarray} 
See the right figure of Figure 6. By Theorem \ref{theorem:main}, 
the relation 
$-\infty < \eta \Leftrightarrow -\infty < \eta = \zeta$ 
holds between QCQP \eqref{eq:QCQP3} and its SDP relaxation \eqref{eq:SDP2} 
with any $\Q \in \SymMat^3$. 
If the objective function ${\scriptsize \begin{pmatrix} \u \\ 1 \end{pmatrix}^T\Q\begin{pmatrix} \u \\ 1 \end{pmatrix}}$ is convex in $\u \in \Real^2$, then \eqref{eq:QCQP3} becomes a convex QCQP, 
so that $\eta = \zeta$ holds. We now consider a nonconvex case. 
For example, consider 
the case where 
\begin{eqnarray*}
& & \Q = \begin{pmatrix} -1 & -1 & 0 \\ -1 & -1 & 0 \\ 0 & 0 & 0 \end{pmatrix}, \
 \begin{pmatrix} \u \\ 1 \end{pmatrix}^T
 \Q
 \begin{pmatrix} \u \\ 1 \end{pmatrix} = -(u_1 + u_2)^2.
\end{eqnarray*}
Obviously, QCQP \eqref{eq:QCQP3} attains the minimum $\zeta = -4$ at 
every boundary points $\u$ of $\BC_\geq$ satisfying $|u_1+u_2| = 2$. 
On the one hand, 
$
\X(t) = {\scriptsize \begin{pmatrix} t & 0 & 1 \\ 0 & t & 1 \\ 1 & 1 & 1 \end{pmatrix} }
$ 
is a feasible solution of SDP \eqref{eq:SDP2} for every $t \geq 2$, and 
$\Q \bullet \X(t) \rightarrow -\infty$  as $t \rightarrow \infty$. 
Therefore, we obtain $\eta = -\infty < \zeta = 4$ in this case. 
\vspace{2mm}\\
{\ } \hspace{2mm}
We note that $\BC_\geq$ is represented by two linear inequalities. We now represent 
the same $\BC_\geq$ by a single quadratic inequality such that 
\begin{eqnarray*}
& & \B^4 = \begin{pmatrix} -1 & -1 & 0 \\ -1 & -1 & 0 \\ 0 & 0 & 4 \end{pmatrix}, \ 
\B^4_{\geq} = \{\u : 0 \leq 4 - (u_1+u_2)^2 \}  = \{\u : -2 \leq u_1+u_2 \leq 2 \}  = \BC_\geq. 
\end{eqnarray*}
In this case, $\eta = \zeta = -4$ holds. 
In general, SDP relaxation of a QCQP depends on the representation of the QCQP. 
In particular, quadratic inequality representation often yields an SDP relaxation that provides a better approximate optimal value (see \cite[Theorem 2.1]{FUJIE1997}). 
\eexamp 


\section{An extension to higher-dimensional instances } 

Several cases of $\BC \subseteq \SymMat^n$ satisfying $\coneJ_+(\BC) \in 
 \wFC(\bGamma^n)$ were provided in 
\cite[Section 4.1]{ARIMA2023} and \cite[Section 6]{ARIMA2024}.  In this section, 
we show how to construct various $\BC$'s  $\subseteq \SymMat^n$ that satisfy Condition (D) 
by combining multiple $\BC$'s $\subseteq \SymMat^3$.  
Recall that we have shown  in Instance 2.6 that 
a combination of two $\BC$'s $\subseteq \SymMat^3$,  
which have the same cardinality and 
 satisfy Condition (D), creates a new 
$\BC \subseteq \SymMat^3$ satisfying Condition (D).  We generalize this technique 
to recursively construct a higher-dimensional $\BC$ satisfying Condition (D). 

Let $1 \leq p \in \Integer$. We assume that both 
$\BC^p = \{\B^{p1},\ldots,\B^{pm}\} \subseteq \SymMat^{n_p}$ and 
$\CC^p = \{\C^{p1},\ldots,\C^{pm}\} \subseteq \SymMat^{\ell_p}$  satisfy Condition (D). 
We note that $(\alpha\B)_\geq = \B_\geq$ for every $\B \in \SymMat^n$ and $\alpha > 0$. 
For simplicity of notation,  we assume that $\alpha_{B^{pk}} = \alpha_{C^{pk}} = 1$ 
$(1 \leq k \leq m)$. 
Thus, our assumption indicates that 
\begin{eqnarray*}
& & \B^{pi} + \B^{pj} \in \SymMat^{n_p}_+ \ (1 \leq i < j \leq m), \ 
\C^{pi}+ \C^{pj} \in \SymMat^{\ell_p}_+ \ (1 \leq i < j \leq m). 
\end{eqnarray*}
Let  $\L^p$ be an $(n_p+ \ell_p)\times n_{p+1}$ matrix. Then, 
we define $\B^{(p+1)i}$ $(1 \leq i \leq m)$ by 
\begin{eqnarray}
\B^{(p+1)i} = (\L^p)^T \begin{pmatrix} \B^{pi} & \O \\ \O^T & \C^{pi} \end{pmatrix} \L^p \in \SymMat^{n_{p+1}}, \label{eq:recursive}
\end{eqnarray}
where $\O$ denotes the $n_p \times \ell_p$ matrix of $0$'s. 
Let $\BC^{p+1} = \{\B^{(p+1)i} : 1 \leq i \leq m \}$. Then, 
\begin{eqnarray*}
\B^{(p+1)i} + \B^{(p+1)j} & = &( \L^p)^T \begin{pmatrix} \B^{pi} & \O \\ \O^T & \C^{pi} \end{pmatrix} \L^p
+ (\L^p)^T \begin{pmatrix} \B^{pj} & \O \\ \O^T & \C^{pj} \end{pmatrix} \L^p \\
& = &( \L^p )^T \begin{pmatrix} \B^{pi}+ \B^{pj}  & \O \\ \O^T & \C^{pi} +\C^{pj} \end{pmatrix} 
\L^p  \in \SymMat^{n_{p+1}}_+ \ ( 1\leq i <  j \leq m ) 
\end{eqnarray*}
holds. Therefore, $\BC^{p+1}$ satisfies Condition (D). Furthermore, by replacing $p$ by $p+1$, 
we can continue this procedure recursively. 
This recursive procedure is highly flexible as it allows 
 arbitrary selection of $\BC^p$ and $\CC^p$ satisfying Condition (D), as well as
any  $(n_p + \ell_p) \times n_{p+1}$ matrix $\L^p$. 
Thus we can create various $\BC$'s in $\SymMat^n$ by this recursive procedure. 

In addition to $\BC$'s $\subseteq \SymMat^3$ described in Section 2, which can be employed 
for the initial $\BC^1$ and $\CC^1$ and also for $\CC^p$ $(p \geq 2)$, we mention some 
$\BC$'s satisfying Condition~(D). 
\begin{description}
\item[(a) ] Every $\BC \subseteq \SymMat^n_+$ satisfies Condition (D). Specifically, every 
$\BC \subseteq  \Real_+$ satisfies Condition~(D).\vspace{-2mm}
\item[(b) ] 
For every $\a \in \Real^{n-1}$ and $\rho > 0$ , let
\begin{eqnarray*}
\B(\a,\rho) = 
\begin{pmatrix} \I  & -\a \\ -\a^T &\a^T\a - \rho^2 \end{pmatrix}. 
\end{eqnarray*}
Then $\B(\a,\rho)_\leq$ forms a ball with the center  $\a$ and radius $\rho$. For every 
$G \subseteq \Integer^{n-1} \times (0,1/2]$ satisfying 
\begin{eqnarray}
\a_1 \not= \a_2 \ 
\mbox{for every pair of distinct $(\a_1,\rho_1)\in G$ and $(\a_2,\rho_2) \in G$}, 
\label{eq:condG}
\end{eqnarray}
define 
$
	\BC(G) = \left\{ \B(\a,\rho) :  (\a,\rho) \in G \right\}. 
$
Then,  we see that 
\begin{eqnarray*}
\lefteqn{ \begin{pmatrix} \I  & -\a_1 \\ -\a_1^T &\a_1^T\a_1 - \rho_1^2 \end{pmatrix} +
\begin{pmatrix} \I  & -\a_2 \\ -\a_2^T &\a_2^T\a_2 - \rho_2^2 \end{pmatrix} }\\
& = & \begin{pmatrix} 2\I  & -(\a_1+\a_2 )\\ -(\a_1+\a_2)^T &\a_1^T\a_1 + \a_2^T\a_2- \rho_1^2 - \rho_2^2\end{pmatrix} \in \SymMat^n_+
\end{eqnarray*}
for every distinct $(\a_1,\rho_1), \ (\a_2,\rho_2) \in G$.  In fact, every leading principal minors of 
the matrix above is nonnegative. (In particular, its determinant can be computed with its 
Schur complement). Therefore, $\BC(G)$ satisfies Condition~(D).\vspace{-2mm}
\item[(c) ] 
Suppose that both $\BC^p$ and $\CC^p$  satisfy Condition (D), but the number $\#\CC^p$ of matrices in 
$\CC^p$ is less than $\#\BC^p$.
To apply \eqref{eq:recursive}, we must  either discard   $\#\BC^p-\#\CC^p$ matrices 
from $\BC^p$ or 
introduce  $\#\BC^p-\#\CC^p$  `dummy'  matrices into $\CC^p$ 
to adjust their cardinalities. 
 If we take a sufficiently large $\lambda \geq 0$, then $\lambda\I$ 
serves as such a dummy matrix that $\CC^p \cup \{ \lambda \I \}$ satisfies Condition (D). 
More precisely, 
let $\lambda_{\rm min}(\B)$ denote 
the minimum eigenvalue of $\B \in \BC$, and $-\lambda \leq \min \{ \lambda_{\rm min}(\B) \
(\B \in \BC), \ 0 \} \leq 0$. Then $\BC \cup \{ \lambda \I\}$ satisfies Condition (D).
\end{description}

\medskip

We now show how to use (b) and (c) for the recursive formula 
\eqref{eq:recursive} to construct $\BC^{p+1}$. Let $\sigma_1 \in \Real$ and $\CC = \{\sigma_1\}$. 
Then $\CC$ satisfies Condition (D) since it consists of a single $1 \times 1$ matrix in $\SymMat^1$. 
For a finite subset $G$ of $\Integer^{n-1} \times (0,1/2]$ 
satisfying \eqref{eq:condG}, let $\BC^p = \BC(G)$, which  satisfies Condition (D). 
Suppose that $\BC^p$ consists of $m \geq 2$ matrices $\B^{p1},\ldots,\B^{pm}$. 
By choosing 
$m-1$ nonnegative numbers $\sigma_i \geq \max\{-\sigma_1,0\}$ $(i=2,\ldots,m)$, 
we expand $\CC = \{\sigma_1\}$ to 
$\CC^p = \{\C^1 = \sigma_1,\C_2 = \sigma_2,\ldots,\C^m = \sigma_m\}$, which satisfies Condition (D)
 by (c) and 
contains the same number of matrices as $\BC^p$. Thus we can apply the recursive formula 
\eqref{eq:recursive} to construct $\BC^{p+1} = \{\B^{(p+1)1},\ldots,\B^{(p+1)m}\}$. 
If we take $\sigma_i = r_i^2 \geq 0$ $(i=1,\ldots,m)$, $n=2$,  $\rho = 1/2$, $G = {G^h}$  and 
$\L^p = {\scriptsize\begin{pmatrix} 0 & 1 & 0 \\ 1 & 0 & 0 \\ 0 & 0 & 1\end{pmatrix}}$, 
then $\BC^{p+1}$ corresponds to $\BC({G^h})$ in Instance~2.4. 

\medskip 

We provide two examples of $\L^p$. For simplicity,  we assume that $n_p = \ell_p = 3$, but 
 the discussion below can be generalized to any 
$n_p, \  \ell_p \geq 3$ in a straightforward manner. 
Let $\lambda_p > 0, \ \mu_p  > 0, \ \lambda_p+\mu_p = 1$. 
For the first example, let
\begin{eqnarray*}
\L^p & = & \begin{pmatrix} \sqrt{\lambda_p}  & 0 & 0 \\
                                       0 & \sqrt{\lambda_p}  & 0 \\
                                       0 & 0 & \sqrt{\lambda_p}  \\
                                       \sqrt{\mu_p}  & 0 & 0\\
                                       0 & \sqrt{\mu_p} & 0 \\
                                       0 & 0 & \sqrt{\mu_p} \end{pmatrix}. 
\end{eqnarray*}
Then, through the linear transformation 
\begin{eqnarray*}
\begin{pmatrix} \u \\ z \end{pmatrix} \in \Real^{2+1} \rightarrow 
\Real^{3+3} \ni \begin{pmatrix} \x^1 \\ \x^2 \end{pmatrix} = \L^p \begin{pmatrix} \u \\ z \end{pmatrix},   
\end{eqnarray*}
the quadratic function ${\scriptsize\begin{pmatrix} \x^1 \\ \x^2 \end{pmatrix}}^T \B^{(p+1)i}{\scriptsize\begin{pmatrix} \x^1 \\ \x^2 \end{pmatrix}}$ in 
${\scriptsize\begin{pmatrix} \x^1 \\ \x^2 \end{pmatrix}}$ is transformed to 
the quadratic function 
\begin{eqnarray*}
\begin{pmatrix} \x^1 \\ \x^2 \end{pmatrix}^T(\L^p)^T\B^{(p+1)i}\L^p\begin{pmatrix} \x^1 \\ \x^2 \end{pmatrix} 
& = & \begin{pmatrix} \u \\ z \end{pmatrix}^T (\lambda_p\B^{pi} + \mu_p\C^{pi})\begin{pmatrix} \u \\ z \end{pmatrix}\\
& = & \lambda_p\begin{pmatrix} \u \\ z \end{pmatrix}^T \B^{pi}\begin{pmatrix} \u \\ z \end{pmatrix} 
+ \mu_p\begin{pmatrix} \u \\ z \end{pmatrix}^T \C^{pi}\begin{pmatrix} \u \\ z \end{pmatrix}
\end{eqnarray*}
in ${\scriptsize\begin{pmatrix} \u \\ z \end{pmatrix}} \in \Real^{2+1}$.
Hence, if we fix $z$ to be  $1$, we obtain a convex combination of the two quadratic functions 
${\scriptsize \begin{pmatrix} \u \\ 1 \end{pmatrix}}^T\B^{pi}{\scriptsize \begin{pmatrix} \u \\ 1 \end{pmatrix}}$ and 
${\scriptsize \begin{pmatrix} \u \\ 1 \end{pmatrix}}^T \C^{pi}{\scriptsize \begin{pmatrix} \u \\ 1 \end{pmatrix}}$. 
Thus this case corresponds to Instance 2.6. 

\bigskip

Now, we consider the case 
\begin{eqnarray*}
\L^p & = & \begin{pmatrix} \sqrt{\lambda_p} & 0 & 0 & 0 & 0 \\
                                       0 & \sqrt{\lambda_p} & 0 & 0 & 0 \\
                                       0 & 0 & 0 & 0 & \sqrt{\lambda_p} \\
                                       0 & 0 & \sqrt{\mu_p} & 0 & 0\\
                                       0 & 0 & 0 & \sqrt{\mu_p} & 0 \\
                                       0 & 0 & 0 & 0 & \sqrt{\mu_p} \end{pmatrix}.                                        
\end{eqnarray*}
In this case, through the linear transformation
\begin{eqnarray*}
\begin{pmatrix} \u^1 \\ \u^2 \\ z \end{pmatrix} \in \Real^{2+2+1} \rightarrow 
\Real^{3+3} \ni \begin{pmatrix} \x^1 \\ \x^2 \end{pmatrix} = \L^p \begin{pmatrix} \u^1 \\ \u^2 \\ z \end{pmatrix},  
\end{eqnarray*}
the quadratic function 
${\scriptsize\begin{pmatrix} \x^1 \\ \x^2 \end{pmatrix}}^T \B^{(p+1)i}
{\scriptsize\begin{pmatrix} \x^1 \\ \x^2 \end{pmatrix}}$ in 
${\scriptsize \begin{pmatrix} \x^1 \\ \x^2 \end{pmatrix}}$ is transformed to 
the quadratic function 
\begin{eqnarray*}
\begin{pmatrix} \x^1 \\ \x^2 \end{pmatrix}^T(\L^p)^T\B^{(p+1)i}\L^p\begin{pmatrix} \x^1 \\ \x^2 \end{pmatrix} 
& = & \begin{pmatrix} \sqrt{\lambda_p} \u^1 \\ \sqrt{\lambda_p} z \\ 
\sqrt{\mu_p} \u^2 \\ \sqrt{\mu_p} z \end{pmatrix}^T 
\begin{pmatrix} \B^{pi}  & \O \\ \O^T & \C^{pi} \end{pmatrix} 
\begin{pmatrix} \sqrt{\lambda_p} \u^1 \\ \sqrt{\lambda_p} z \\ 
\sqrt{\mu_p} \u^2 \\ \sqrt{\mu_p} z \end{pmatrix} \\
& = & \lambda_p\begin{pmatrix} \u^1 \\ z \end{pmatrix}^T \B^{pi}\begin{pmatrix} \u^1 \\ z \end{pmatrix} 
+ \mu_p\begin{pmatrix} \u^2 \\ z \end{pmatrix}^T \C^{pi}\begin{pmatrix} \u^2 \\ z \end{pmatrix}
\end{eqnarray*}
in ${\scriptsize\begin{pmatrix} \u^1 \\ \u^2 \\ z \end{pmatrix}} \in \Real^{2+2+1}$. 
Therefore, if we fix  $z$ to be  $1$,
we obtain a convex combination of the quadratic function 
${\scriptsize\begin{pmatrix} \u^1 \\ 1 \end{pmatrix}}^T \B^{pi}
{\scriptsize\begin{pmatrix} \u^1 \\ 1 \end{pmatrix}}$ in $\u^1 \in \Real^2$ and 
the quadratic function 
${\scriptsize\begin{pmatrix} \u^2 \\ 1 \end{pmatrix}}^T \C^{pi}
{\scriptsize\begin{pmatrix} \u^2 \\ 1 \end{pmatrix}}$ in 
$\u^2 \in \Real^2$.  

\rema
To add a linear equality $\A\u = \b$ to QCQP \eqref{eq:QCQP3}, we introduce 
$\B$ $=$ $- (\A, -\b)^T$ $(\A,-\b)$ $\in$ $\SymMat^n$, where $\A$ denotes an $(n-1)\times \ell$ matrix, 
and $\b \in \Real^\ell$. We see that 
\begin{eqnarray*}
\B_\geq  = \left\{ \u \in \Real^{n-1} : \begin{pmatrix} \u \\ 1 \end{pmatrix}^T
\B\begin{pmatrix} \u \\ 1 \end{pmatrix} \geq 0 \right\} 
=  \left\{ \u \in \Real^{n-1} : -\parallel \A\u - \b \parallel^2  \geq 0 \right\}.
\end{eqnarray*}
Hence $\A \u = \b$ if and only if $\u \in \B_\geq$. It is known that 
if $\BC$ satisfies Condition (B), then so does $\BC' = \BC \cup \{ \B \}$. 
See \cite[Section 4.4]{ARIMA2023} for more details. 
\erema


\section{Computing  a QCQP optimal solution from an optimal solution of its SDP relaxation} 

Throughout this section, we assume \vspace{-2mm}
\begin{description}
\item[(a) ]  $\BC = \{\B^1,\ldots,\B^m\} \subseteq \SymMat^n$ satisfies Condition (B),\vspace{-2mm} 
\item[(b) ]  The SDP relaxation~\eqref{eq:SDP2} of QCQP~\eqref{eq:QCQP2} has an optimal solution 
$\overline{\X} \in \SymMat^n$,\vspace{-2mm}
\item[(c) ] $\X = \overline{\X} \in \SymMat^n$ satisfies
 {\em the KKT (Karush-Kuhn-Tucker) stationary condition}: 
there exists a $(\bar{t},\bar{\y},\overline{\Y}) \in \Real \times \Real^m \times \SymMat^n$ such that 
\begin{eqnarray}
\left. 
\begin{array}{l}
\X \in \SymMat^n_+, \ \B^k \bullet \X \geq 0 \ (1 \leq k \leq m), \ \H \bullet \X = 1 \ \mbox{(primal feasibility)}, \\[3pt]
 \bar{\y} \geq \0, \ \Q - \H \bar{t} - \sum_{k=1}^m \bar{y}_k \B^k = \overline{\Y} \in \SymMat^n_+ 
\ \mbox{(dual feasibility)}, \\[3pt]
 \bar{y}_k (\B^k \bullet \X)= 0 \ (1 \leq k \leq m), \ 
\overline{\Y} \bullet \X = 0 \ \mbox{(complementarity)}. 
\end{array}
\right\}  
\label{eq:KKT}
\end{eqnarray}
\vspace{-2mm}
\end{description}
In (c), $(\bar{t},\bar{\y},\overline{\Y}) \in \Real \times \Real^m\times \SymMat^n$ corresponds 
 to an optimal solution of the dual of SDP~\eqref{eq:SDP2}. We note that (c) $\Rightarrow$ (b). 
If $\H \in \SymMat^n_+$, 
in particular,  if $\H = \mbox{diag}(0,\ldots,0,1) \in \SymMat^n_+$ as in QCQP~\eqref{eq:QCQP3}, 
then (b) $\Rightarrow$ (c) \cite[Theorem 2.1]{KIM2022}. 
By Theorem~\ref{theorem:main}, the optimal values of SDP~\eqref{eq:SDP2} and QCQP~\eqref{eq:QCQP2} coincide, {\it i.e.}, 
$\eta = \zeta$. 

We will describe a numerical method for computing an optimal solution  
$\widetilde{\X} \in \bGamma^n$
of QCQP~\eqref{eq:QCQP2} from $\overline{\X} \in \SymMat^n_+$. 
Recall that all $2$-dimensional QCQP  instances in Section 2.3 
satisfy Condition (D) and that a recursive procedure is provided 
for constructing higher-dimensional QCQP instances satisfying 
Condition (D)  in Section 3. Since Condition (D) implies Condition (B) by 
Theorem~\ref{theorem:main} (iv), we can apply 
the method to those instances. 

\medskip

Let $K_0 = \{ k  : \B^k \bullet \overline{\X} = 0\}$. Then we have either \vspace{-1mm}
\begin{description}
\item{(i) } $k \in K_0 \not= \emptyset$ for some $k \in \{1,\ldots,m\}$.\vspace{-2mm}
\item{(ii) } $K_0 = \emptyset$.
\end{description}
We first deal with case (i). In this case, the method is based on the following lemma and its 
constructive proof. 

\lemm \label{lemma:YeZhang}
(\cite[Lemma 2.2]{YE2003}, see also \cite[Proposition 3]{STURM2003}) 
Let $\B \in \SymMat^n$ and $\overline{\X}\in \SymMat^n_+$ with 
rank$\overline{\X} = r$. Suppose that $\B \bullet\overline{\X} \geq 0$. Then, there exists a rank-1 
decomposition of $\overline{\X}$ such that $\overline{\X} = \sum_{i=1}^r\x_i\x_i^T$ 
and $\B \bullet \x_i\x_i^T \geq 0$ 
$(1 \leq i \leq r)$. If, in particular, $\B \bullet\overline{\X}  = 0$, then 
$\B\bullet\x_i\x_i^T = 0$ $(1 \leq i \leq r)$. 
\elemm

By Lemma~\ref{lemma:YeZhang}, 
there exists a rank-1 decomposition of $\overline{\X}$ such that 
$\overline{\X} = \sum_{i=1}^r\x_i\x_i^T$ 
and $\B^k \bullet x_i x_i^T= 0$ 
({\it i.e.}, $\x_i\x_i^T \in \coneJ_0(\B^k)$) $(1 \leq i \leq r)$. 
By assumption (a), $\x_i\x_i^T \in \coneJ_+(\BC) \ (1 \leq i \leq r)$.
Since 
$1 = \H \bullet \overline{\X} = \sum_{i=1}^r \H \bullet \x_i\x_i^T$, there exist a 
$\tau \geq 1/r$ and a $j \in \{1,\ldots,r\}$ such that $\H \bullet \x_j\x_j^T= \tau$. 
Let $\widetilde{\X} = \x_j\x_j^T/\tau$. Since $ \x_j\x_j^T\in \coneJ_+(\BC)$, 
$\overline{\X} \in \coneJ_+(\BC)$. We also see 
that $\H \bullet \widetilde{\X} = \H \bullet  (\x_j^T\x_j/\tau )= 1$. Hence $ \widetilde{\X} \in \bGamma^n$ 
is a rank-$1$ feasible solution of SDP~\eqref{eq:SDP2}. Furthermore, 
we see from $\overline{\Y} \in \SymMat^n_+$ and $\x_i\x_i^T \in \SymMat^n_+$ 
$(1 \leq i \leq r)$ that 
\begin{eqnarray*}
0 \leq \overline{\Y}\bullet \widetilde{\X} = \frac{\overline{\Y}\bullet \x_j\x_j^T}{\tau} \leq 
\frac{\overline{\Y} \bullet \sum_{i=1}^m\x_i\x_i^T}{\tau} = 
\frac{\overline{\Y} \bullet \overline{\X}}{\tau} = 0.  
\end{eqnarray*}
Hence, $\widetilde{\X} \in \bGamma^n$ is a rank-1 optimal solution of SDP~\eqref{eq:SDP2}  
and it is an optimal solution of  QCQP~\eqref{eq:QCQP2} with the 
same objective value $\Q \bullet \overline{\X}$. 

We now describe how to compute the rank-$1$ decomposition of 
$\overline{\X}$ such that $\overline{\X} = \sum_{i=1}^r\x_i\x_i^T$ and 
$\B\bullet\x_i\x_i^T = 0$ $(1 \leq i \leq r)$ based on the constructive proof of 
Lemma~\ref{lemma:YeZhang}  in \cite{STURM2003,YE2003}. 
Let $\overline{\X} = \sum_{i=1}^r\x_i\x_i^T$  be an arbitrary rank-$1$ decomposition. 
If $\B\bullet\x_i\x_i^T = 0$ $(1 \leq i \leq r)$, then we are done. Otherwise, 
there exist $i$ and $j$ such that $\B\bullet\x_i\x_i^T < 0 <  \B\bullet\x_j\x_j^T$, say $i=1$ and $j=2$. 
We consider the following quadratic equation in $s$:
\begin{eqnarray*}
0 & = & \B \bullet (s\x_1+\x_2) (s\x_1+\x_2)^T = (\B\bullet\x_1\x_1^T)s^2 + 2 (\B\bullet\x_1\x_2^T) s + 
\B\bullet\x_2\x_2^T. 
\end{eqnarray*}
Since $\B\bullet\x_1\x_1^T \times \B\bullet\x_2\x_2^T < 0$, this equation has two distinct roots with 
opposite signs. Let $\bar{s}$ be one of the roots. Let
\begin{eqnarray*}
& & \bar{\x}_1 = \frac{\bar{s}}{\sqrt{\bar{s}^2+1}}\x_1 + \frac{1}{\sqrt{\bar{s}^2+1}}\x_2 \ \mbox{and }
 \bar{\x}_2 = -\frac{1}{\sqrt{\bar{s}^2+1}}\x_1 + \frac{\bar{s}}{\sqrt{\bar{s}^2+1}}\x_2.
\end{eqnarray*}
Then, we have 
\begin{eqnarray*}
& &  \B\bullet \bar{\x}_1\bar{\x}_1^T  =  0, \ 
 \bar{\x}_1\bar{\x}_1^T  + \bar{\x}_2\bar{\x}_2^T 
 =  \x_1\x_1^T + \x_2\x_2^T.  
\end{eqnarray*}
Now replace $\x_1$ with $\bar{\x}_1$ and $\x_2$ with $\bar{\x}_2$. Then we still have  the 
rank-$1$ decomposition $\overline{\X} = \sum_{i=1}^r\x_i\x_i^T$. If $\B \bullet \x_i\x_i^T \not= 0$ 
for some $i \geq 2$,  we continue this procedure recursively till all 
$\B \bullet \x_i\x_i^T = 0$ $(0 \leq i \leq r)$ are attained. 
To compute an optimal solution 
$\widetilde{\X} \in \bGamma^n$ of QCQP~\eqref{eq:QCQP2}, we can terminate the procedure once we 
find $\x_j \in \Real^n$  such that $\B\bullet\x_j\x_j^T = 0$ and $\tau \equiv \H\bullet\x_j\x_j^T \geq 1/r$. 
In this case, 
 $\widetilde{\X} = \x_j\x_j^T/\tau$ is an optimal solution of QCQP~\eqref{eq:QCQP2}. 

\medskip

We now consider case (ii). 
 In the KKT condition~\eqref{eq:KKT}, $\X  = \overline{\X}$ satisfies 
$\B^k \bullet  \overline{\X}  > 0$ and $\bar{y}_k = 0$ $(1 \leq k \leq m)$. This implies that 
$\overline{\X}$ is an optimal solution of a simple SDP with the single equality constraint 
$ \H \bullet \X = 1$: 
\begin{eqnarray}
\tilde{\eta} = \inf \left\{ \Q \bullet \X : \X \in \SymMat^n_+, \ \H \bullet \X = 1 \right\}.   \label{eq:SDP4} 
\end{eqnarray}
The KKT stationary condition for SDP~\eqref{eq:SDP4}  is written as 
\begin{eqnarray}
\left. 
\begin{array}{l}
\X \in \SymMat^n_+, \ \H \bullet \X = 1 \ \mbox{(primal feasibility)} \\[3pt]
 \Q - \H \bar{t} = \overline{\Y} \in \SymMat^n_+ 
\ \mbox{(dual feasibility)} \\[3pt]
\overline{\Y} \bullet \X = 0 \ \mbox{(complementarity)} 
\end{array}
\right\},
\label{eq:KKT2}
\end{eqnarray}
for some $(\tilde{t},\widetilde{\Y}) \in \Real \times \SymMat^n$,
which serves as a sufficient condition for 
$\X \in \SymMat^n$ to be an optimal solution of SDP~\eqref{eq:SDP4}. 
Now we compute a rank-$1$ optimal solution of SDP~\eqref{eq:SDP4} from $\overline{\X}$. 
Let $r = \mbox{rank}\overline{\X}$. If $r = 1$, then 
we have done. So assume $r \geq 2$. Let $\overline{\X} = \sum_{i=1}^r \x_i\x_i^T$ be a 
rank-$1$ decomposition of $\overline{\X}$. 
It follows from~\eqref{eq:KKT} with $\X = \overline{\X}$ that 
\begin{eqnarray*}
0 = \overline{\Y} \bullet \overline{\X} =  \overline{\Y} \bullet (\sum_{i=1}^r \x_i\x_i^T) = \sum_{i=1}^r (\overline{\Y}  \bullet \x_i\x_i^T). 
\end{eqnarray*}
Since $\overline{\Y} \in \SymMat^n_+$ and $\x_i\x_i^T \in \SymMat^n_+$, 
$\overline{\Y} \bullet \x_i\x_i^T \geq 0$ $(1 \leq i \leq r)$. Hence the identity above implies 
$\overline{\Y} \bullet \x_i\x_i^T = 0$ $(1 \leq i \leq r)$. We also see from 
$1 = \H \bullet \overline{\X} = \sum_{i=1}^r (\H \bullet \x_i\x_i^T)$ that there exists a $j$ such that 
$\tau \equiv \H \bullet \x_j\x_j^T \geq  1/r$. Define $\widetilde{\X} = \x_j\x_j^T/\tau$. Then 
$\X =  \widetilde{\X}$ satisfies~\eqref{eq:KKT2}. Thus $\widetilde{\X}$ is a rank-$1$ optimal solution of 
SDP~\eqref{eq:SDP4}. If $\B^k \bullet \widetilde{\X}\geq 0$ $(1 \leq k \leq m)$ then 
$\widetilde{\X}$ is a rank-$1$ optimal solution of 
SDP~\eqref{eq:SDP2}, hence an optimal solution of QCQP~\eqref{eq:QCQP2}.  
Otherwise, we can take an optimal solution $\widehat{\X}$ of 
SDP~\eqref{eq:SDP2} such that $\{ k : \B^k \bullet \widehat{\X} = 0 \} \not= \emptyset$ 
as a convex combination of $\overline{\X}$ and 
$\widetilde{\X}$, which leads to case (i). 

\begin{figure}[t!]   \vspace{-30mm} 
\begin{center}
\includegraphics[height=150mm]{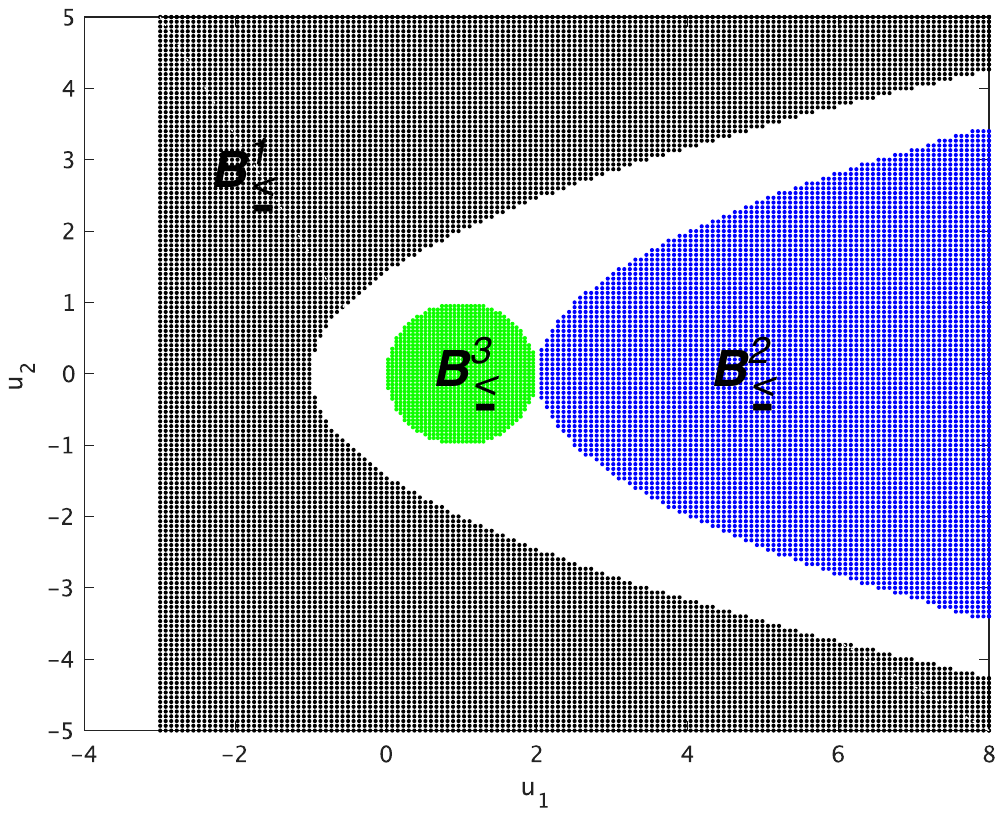}
\end{center}

\vspace{-35mm} 
\caption{
The unshaded region  represents the feasible region $\BC_\geq$ of QCQP~\eqref{eq:QCQP4}. 
}
\end{figure}

\smallskip

To illustrate how the method works, we consider the following instance.
\examp
We simultaneously describe 6 QCQP problems with quadratic objective functions $q^k: \Real^2 \rightarrow \Real$ 
$(1 \leq k \leq 6)$ over a common quadratic inequality feasible region: 
\begin{eqnarray}
\eta & = & \inf\left\{q^k(\u) : \u \in \Real^2, \ -2 \leq 2u_1 -u_2^2 <= 4, \ (u_1-1)^2+u_2^2 \geq 1 \right\} 
\nonumber  \\
       & = & 
 \inf \left\{ \begin{pmatrix}\u \\ 1\end{pmatrix}^T\Q^k\begin{pmatrix}\u \\ 1\end{pmatrix} : 
 \u \in \BC_\geq \right\} \label{eq:QCQP4}  \\ 
       & = & \inf \left\{\Q^k\bullet\X : \X\in\bGamma^3\cap\coneJ_+(\BC), \ 
                  \mbox{diag}(0,0,1)\bullet\X = 1 \right\}, \label{eq:QCQP5}
\end{eqnarray}
where
\begin{eqnarray*}
& & \B^1 = {\scriptsize\begin{pmatrix}  0  & 0  &1 \\ 0 &  -1 & 0 \\ 1  & 0 & 2 \end{pmatrix}}, \ 
 \B^2 = {\scriptsize\begin{pmatrix}  0  & 0  &-1 \\ 0 &  +1 & 0 \\ -1  & 0 & 4 \end{pmatrix}}, \ 
 \B^3 = {\scriptsize\begin{pmatrix}  1  & 0  &-1 \\ 0 &  +1 & 0 \\ -1  & 0 & 0 \end{pmatrix}}, \ 
 \BC = \{\B^1,\B^2,\B^3\}, \\
 & & q^1(\u) = (u_1-2)^2 + (u_2-1)^2 
 = {\scriptsize \begin{pmatrix}\u \\ 1 \end{pmatrix}}^T\Q^1{\scriptsize \begin{pmatrix}\u \\ 1 \end{pmatrix}}, \ \Q^1 =  
 {\scriptsize \begin{pmatrix} 1&0&-2\\0&1&-1\\-2&-1&5\end{pmatrix}}, \\
  & & q^2(\u) = (u_1+3)^2 + u_2^2 = {\scriptsize \begin{pmatrix}\u \\ 1 \end{pmatrix}}^T\Q^2{\scriptsize \begin{pmatrix}\u \\ 1 \end{pmatrix}}, \ \Q^2 =  
 {\scriptsize \begin{pmatrix} 1&0&3\\0&1&0\\3&0&9\end{pmatrix}}, \\
  \end{eqnarray*}
 \begin{eqnarray*} 
  & & q^3(\u) = 2u_1= {\scriptsize \begin{pmatrix}\u \\ 1 \end{pmatrix}}^T\Q^3{\scriptsize \begin{pmatrix}\u \\ 1 \end{pmatrix}}, \ \Q^3 =  
 {\scriptsize \begin{pmatrix} 0&0&1\\0&0&0\\1&0&0\end{pmatrix}}, \\
  & & q^4(\u) = 0= {\scriptsize \begin{pmatrix}\u \\ 1 \end{pmatrix}}^T\Q^4{\scriptsize \begin{pmatrix}\u \\ 1 \end{pmatrix}}, \ \Q^4 =  \O, \\
  & & q^5(\u) = (u_1+4u_2-4)^2 = {\scriptsize \begin{pmatrix}\u \\ 1 \end{pmatrix}}^T\Q^5{\scriptsize \begin{pmatrix}\u \\ 1 \end{pmatrix}}, \ \Q^5 =  
 {\scriptsize \begin{pmatrix} 1&4 &-4\\4&16&-16\\-4&-16&16\end{pmatrix}} \\
& & q^6(\u) = (u_1-3)^2= {\scriptsize \begin{pmatrix}\u \\ 1 \end{pmatrix}}^T\Q^6{\scriptsize \begin{pmatrix}\u \\ 1 \end{pmatrix}}, \ \Q^6 =  
 {\scriptsize \begin{pmatrix} 1 &0&-3\\0&0&0\\-3&0&9\end{pmatrix}}. 
 \end{eqnarray*}
See Figure 13 for the feasible region $\BC_\geq$ of QCQP~\eqref{eq:QCQP4}.

\smallskip

Table 1 presents a summary of the numerical results obtained for solving QCQP~\eqref{eq:QCQP4}. 
The optimal solutions of QCQP~\eqref{eq:QCQP4} can be easily obtained from Figure 13 
for $k=1,\ldots,6$, 
as shown in Table 1. The 
optimal solutions $\overline{\X}$ of the SDP relaxation were computed by SeDuMi \cite{STURM99}. 
For $k=1$ or $k=2$, QCQP~\eqref{eq:QCQP4} has a unique optimal solution, and 
a rank-$1$ solution $\overline{\X}$ was obtained by just solving the SDP relaxation. 
For $k=3$, 
QCQP~\eqref{eq:QCQP4} has a unique optimal solution $(-1,0)$, but 
the computed optimal solution $\overline{\X}$ of the SDP relaxation has rank-$2$ and case (i) occurred. 
For $k=4$,  $5$ or $6$,  the SDP relaxation 
as well as QCQP~\eqref{eq:QCQP4} have 
multiple optimal solutions, 
and rank$\overline{\X} \geq 2$.
For $k=4$ or $k=5$, 
the rank-$1$ decomposition 
$\overline{\X} = \sum_{i=1}^3\x_i\x_i^T$ in case (ii) 
led to an optimal solution $\widetilde{\X} = \x_j\x_j^T/\tau$ of QCQP~\eqref{eq:QCQP5} for some $j$ 
and $\tau = \mbox{diag}(0,0,1)\bullet \x_j\x_j^T \geq 1/3$. 
For $k=6$, 
case (ii) occurred at the optimal solution $\overline{\X}$ of the SDP relaxation, but its rank-$1$ decomposition 
$\overline{\X} = \sum_{i=1}^3\x_i\x_i^T$ yielded no feasible rank-$1$ solution $\widetilde{\X}$ and 
the recursive procedure of case (i) was carried out from a convex combination of $\x_j\x_j^T/(\mbox{diag}(0,0,1)\bullet\x_j\x_j^T)$  for some $j$ and $\overline{\X}$.

\begin{table}[h!]
\begin{center}
\begin{tabular}{|c|c|c|c|c|c|c|c|c|}
\hline
               &  \multicolumn{2}{|c|}{QCQP} & \multicolumn{6}{|c|}{SDP, $\overline{\X}$: Computed Optimal Solution }  \\
$k$  &  Opt.Sol. & Opt.Val. & Opt.Val & Rank & $\B^1\bullet\overline{\X}$ & $\B^2\bullet\overline{\X}$ & 
$\B^3\bullet\overline{\X}$ & Case  \\
\hline 
1              & $\u=(2,1)$                & 0     & $0.00$      & 1    & 5.00   & 1.00 & 1.00 & (ii) \\
\hline
2              & $\u=(-1,0)$               & 4     & $4.00$       & 1   & $0.00$ & 6.00 & 3.00 & (i) \\
\hline 
3              & $\u=(-1,0)$               & -2    & $-2.00$     & 2   & $0.00$ & 6.00 &  5.26 & (i) \\
\hline 
4             &   $\forall \u \in S^4$   & 0    & $0.00$                           & 3   & 2.00  & 3.99 & 2.28  & (ii) \\
\hline 
5              &  $\forall \u \in S^5$   & 0    & $0.00$                      &  2 & 1.68   & $4.32$ & 1.42 & (ii) \\
\hline 
6              &   $\forall \u \in S^6$   & 0    & $0.00$  &  2    & 3.45 & 2.55 & 7.55 & (ii) $\rightarrow$ (i) \\
\hline                
\end{tabular}	
\end{center}
\caption{ 
$S^4 = \BC_\geq$. 
$S^5 = \left\{ \u \in \BC_\geq:  u_1+4u_2=4 \right\}$. 
 $S^6 = \left\{ \u \in \BC_\geq:  u_1 = 3 \right\}$. 
 }
\end{table}
\eexamp


\section{Concluding remarks}

If $\BC \subseteq \SymMat^n$ satisfies $\coneJ_+(\BC) \in 
 \wFC(\bGamma^n)$, then QCQP \eqref{eq:QCQP2} is equivalent to 
its SDP relaxation~\eqref{eq:SDP2} whose optimal value and solution can be easily   computed.
In Section 3, we have shown a recursive procedure for constructing 
various $\BC$'s  $\subseteq \SymMat^n$ satisfying Condition (D). 
The class of QCQPs with such $\BC$'s appears to be quite broad, 
at least in theory.  It should be noted 
that Condition (D) is merely sufficient for $\coneJ_+(\BC) \in 
 \wFC(\bGamma^n)$. 
For example, we could relax Condition (D) to 
 \begin{description}
 \item[Condition (D)': ] If $\A, \ \B \in \BC$ and $\A \not= \B$ then $\alpha\A+\beta\B \in \SymMat^n_+$ 
for some nonzero $(\alpha,\beta) \in \Real^2$ (\cite[Proposition 1]{ARGUE2023}). 
\end{description}
There is still a gap, however,  between Condition (D)' and  the condition $\coneJ_+(\BC) \in 
 \wFC(\bGamma^n)$ 
(\cite[Theorem 1]{ARGUE2023} and \cite[Example 6.1]{ARIMA2024}).  
See \cite[Section 3]{ARIMA2024} 
where various sufficient conditions for 
$\coneJ_+(\BC) \in  \wFC(\bGamma^n)$ and 
their relationships are discussed, under moderate assumptions including the 
Slater constraint qualification  that 
$\coneJ_+(\BC) \subseteq \SymMat^n_+$ contains a positive definite matrix.



\begin{thebibliography}{10}

\bibitem{ARGUE2023}
C.~J. Argue, F.~{Kilin{\c c}-Karzan}, and A.L. Wang.
\newblock Necessary and sufficient conditions for rank-one-generated cones.
\newblock {\em Math. Oper. Res.}, 48(1):100--126, 2023.

\bibitem{ARIMA2023}
N.~Arima, S.~Kim, and M.~Kojima.
\newblock Further development in convex conic reformulation of geometric
  nonconvex conic optimization problems.
\newblock {\em SIAM J. Optim.}, 34(4):3194--3211, August 2024.

\bibitem{ARIMA2024}
N.~Arima, S.~Kim, and M.~Kojima.
\newblock Exact {SDP} relaxations for a class of quadratic programs with finite
  and infinite quadratic constraints.
\newblock Technical Report arXiv:2409.07213, September 2024.

\bibitem{AZUMA2022}
G.~Azuma, Fukuda M., S.~Kim, and M.~Yamashita.
\newblock Exact {SDP} relaxations for quadratic programs with bipartite graph
  structures.
\newblock {\em J. of Global Optim.}, 86:671--691, 2023.

\bibitem{FUJIE1997}
T.~Fujie and M.~Kojima.
\newblock Semidefinite programming relaxation for nonconvex quadratic programs.
\newblock {\em J. of Global Optim.}, 10:367--368, 1997.

\bibitem{HILDEBRAND2016}
R.~Hildebrand.
\newblock Spectrahedral cones generated by rank 1 matrices.
\newblock {\em J. Global Optim}, 64:349--397, 2016.

\bibitem{KIM2003}
S.~Kim and M.~Kojima.
\newblock Exact solutions of some nonconvex quadratic optimization problems via
  {SDP} and {SOCP} relaxations.
\newblock {\em Comput. Optim. Appl.}, 26(2):143--154, 2003.

\bibitem{KIM2022}
S.~Kim and M.~Kojima.
\newblock Strong duality of a conic optimization problem with a single
  hyperplane and two cone constraints strong duality of a conic optimization
  problem with a single hyperplane and two cone constraints.
\newblock {\em Optimization}, 74(1):33--53, 2025.

\bibitem{KIM2020}
S.~Kim, M.~Kojima, and K.~C. Toh.
\newblock A geometrical analysis of a class of nonconvex conic programs for
  convex conic reformulations of quadratic and polynomial optimization
  problems.
\newblock {\em SIAM J. Optim.}, 30:1251--1273, 2020.

\bibitem{MURTY87}
K.~G. Murty and S.~N. Kabadi.
\newblock Some \mbox{NP-complete} problems in quadratic and non-linear
  programming.
\newblock {\em Math. {P}rogram.}, 39:117--129, 1987.

\bibitem{SHOR1987}
N.~Z. Shor.
\newblock Quadratic optimization problems.
\newblock {\em Soviet Journal of Computer and Systems Sciences}, 25:1--11,
  1987.

\bibitem{Shor1990}
N.~Z. Shor.
\newblock Dual quadratic estimates in polynomial and boolean programming.
\newblock {\em Ann. Oper. Res.}, 25:163--168, 1990.

\bibitem{SOJOUDI2014}
S.~Sojoudi and J.~Lavaei.
\newblock Exactness of semidefinite relaxations for nonlinear optimization
  problems with underlying graph structure.
\newblock {\em SIAM J. Optim.}, 24(4):1746--1778, 2014.

\bibitem{STURM99}
J.~F. Sturm.
\newblock {S}e{D}u{M}i 1.02, a {MATLAB} toolbox for optimization over symmetric
  cones.
\newblock {\em {O}ptim. {M}ethods and {S}oftw.}, 11{\&}12:625--653, 1999.

\bibitem{STURM2003}
J.~F. Sturm and S.~Zhang.
\newblock On cones of nonnegative quadratic functions.
\newblock {\em Math. Oper. Res.}, 28(2):246--267, 2003.

\bibitem{YE2003}
Y.~Ye and S.~Zhang.
\newblock New results on quadratic minimization.
\newblock {\em SIAM J. Optim.}, 14:245--267, 2003.

\bibitem{ZHANG2000}
S.~Zhang.
\newblock Quadratic optimization and semidefinite relaxation.
\newblock {\em Math. {P}rogram.}, 87:453--465, 2000.

\end{thebibliography}

\end{document}